\documentclass[12pt,reqno]{amsart}
\usepackage{amssymb,amsfonts,amsthm}
\usepackage{ulem}
\usepackage{amsthm}
\usepackage{amssymb}
\usepackage{mathrsfs}
\usepackage{frcursive}
\usepackage[english,american]{babel}
\usepackage{enumerate}
\usepackage{graphicx}
\usepackage{tikz} 
\usepackage{setspace}
\usepackage{hyperref}
\usepackage{color}

\definecolor{Red}{cmyk}{0,1,1,0.2}

\usepackage{amsfonts,amsmath,layout}

\newcommand{\R}{\mathbb R}

\def\R{\mathbb R}

\def\E{\mathbb E}
\def\P{\mathbb P}

\newcommand{\be}{\begin{equation}}
\newcommand{\ee}{\end{equation}}

\def\1{{\bf 1}}

\def\ds{\displaystyle}

\newtheorem{Theorem}{Theorem}[section]
\newtheorem{Definition}[Theorem]{Definition}
\newtheorem{Proposition}[Theorem]{Proposition}

\begin{document}
\title[HJB equations with local-time Kirchhoff's boundary condition]{Comparison principle for Walsh's spider HJB equations with non linear local-time Kirchhoff's boundary condition}
\space
\space
\author[Isaac Ohavi]{Isaac Ohavi $^\aleph$
\\
$^\aleph$Hebrew University of Jerusalem, Department of  Mathematics and Statistics, Israël}
\address{ Hebrew University of Jerusalem, Department of Mathematics and Statistics, Israël}
\email{isaac.ohavi@mail.huji.ac.il \& isaac.ohavi@gmail.com}
\thanks{This research  was supported by the GIF grant 1489-304.6/2019.\\
I am grateful to Miguel Martinez$^{\beth}$ for the several discussions and collaborations on spider diffusions and their related partial differential equations. I am also grateful to Ioannis Karatzas $^\gimel$, with the various exchanges we were able to have with the different members of the seminary: Probability and Control theory - Columbia University, New York.\\
$^\beth$ Université Gustave Eiffel, LAMA UMR 8050, France. E-mail: miguel.martinez@univ-eiffel.fr\\
$^\gimel$ Department of Mathematics, Columbia University, New York, NY 10027, USA. E-mail: ik1@columbia.edu.}
\dedicatory{Version: \today}
\maketitle
\begin{abstract} The main purpose of this work is to obtain a comparison principle for viscosity solutions of a system of elliptic Walsh's spider Hamilton-Jacobi-Bellman (HJB) equations, possessing a new boundary condition called {\it non linear local-time Kirchhoff's boundary condition}. The main idea is to build test functions in a neighborhood of the vertex with well-designed coefficients. The key point is to impose a 'local-time' derivative at the vertex absorbing the error term induced by - what we decide to call here - the {\it Kirchhoff's speed of the Hamiltonians}. As a consequence, we obtain a comparison theorem for HJB systems posed on star-shaped networks, with non linear Kirchhoff's boundary condition and non vanishing viscosity at the vertex.
\end{abstract}
\section{Introduction}
We are given an integer number $I$ (with $I\ge 2$) and star-shaped compact network: $$\displaystyle \mathcal{N}_R=\bigcup_{i=1}^I\mathcal{R}_i,$$ that consists of $I$ compact rays $\mathcal{R}_i\cong [0,R]$ ($R>0$) emanating from a junction point $\bf 0$.

The main target of this work is to obtain a comparison theorem (thus uniqueness) for continuous viscosity solution of the following Walsh's spider Hamilton-Jacobi-Bellman system - $\mathcal{W}_{alsh}(\mathbb{S})$ - having a new boundary condition at the vertex $\bf 0$, called {\it non linear local-time Kirchhoff's boundary condition:}
\label{sec intro}
\begin{eqnarray}\label{eq PDE 0}
\nonumber &\mathcal{W}_{alsh}(\mathbb{S}):=\\
&\begin{cases}
\textbf{HJB equations parameterized by the local-time on each ray:}\\ 
\lambda u_i(x,l)+\underset{\beta_i\in \mathcal{B}_i}{\sup}\Big\{-\sigma_i(x,l,\beta_i)\partial^2_xu_i(x,l)+
b_i(x,l,\beta_i)\partial_xu_i(x,l)+\\
h_i(x,l,\beta_i)\Big\}=0,~~(x,l)\in(0,R)\times (0,K),\\
\textbf{Non linear local-time Kirchhoff's boundary condition at } \bf 0:\\
\partial_lu(0,l)+\underset{ \vartheta \in \mathcal{O}}{\inf} \Big\{\ds\sum_{i=1}^I\mathbb{S}_i(l,\vartheta)\partial_xu_i(0,l)+h_0(l,\vartheta)\Big\}=0,~~l\in(0,K),\\
\textbf{Dirichlet boundary conditions outside } \bf 0:\\ 
u_i(R,l)=\chi_i(l),~~l\in[0,K],\\
u_i(x,K)=\mathcal{T}_i(x),~~x\in[0,R],\\
\textbf{Continuity condition at } \bf 0:\\ 
\forall (i,j)\in[\![1,I]\!]^2,~~\forall l\in [0,K],~~u_i(0,l)=u_j(0,l).
\end{cases}
\end{eqnarray}
In all of this work, $\mathcal W_{alsh}(\mathbb S)$ will be used as a label to refer to system \eqref{eq PDE 0}.
In order to simplify our study, we have assumed in our framework that all the rays $\mathcal{R}_i=[0,R]\times \{i\},~i\in[\![1,I]\!]$, have the same length $R>0$, and that Dirichlet boundary conditions $\chi_i$ hold at $x=R$ and $\mathcal{T}_i$ at $l=K$. A more general setting could be treated with similar tools: one could for instance consider more general rays, and/or a mix of local-time Kirchhoff's and Dirichlet boundary conditions at $x=R$, $l=K$, etc. 

System \eqref{eq PDE 0} - $\mathcal{W}_{alsh}(\mathbb{S})$ - can be interpreted as the non linear version in the elliptic framework, of the following linear parabolic system:
\begin{eqnarray}\label{eq : pde with l}
\begin{cases}
\textbf{Linear parabolic equation parameterized}\\
\textbf{by the local-time on each ray:}\\
\partial_tu_i(t,x,l)-\sigma_i(t,x,l)\partial_x^2u_i(t,x,l)
+b_i(t,x,l)\partial_xu_i(t,x,l)\\
+c_i(t,x,l)u_i(t,x,l)=f_i(t,x,l),~~(t,x,l)\in (0,T)\times (0,R)\times(0,K),\\
\textbf{Linear local-time Kirchhoff's boundary condition at }\bf 0:\\
\partial_lu(t,0,l)+\displaystyle \sum_{i=1}^I \alpha_i(t,l)\partial_xu_i(t,0,l)=\phi(t,l),~~(t,l)\in(0,T)\times(0,K),\\
\medskip
\textbf{Dirichlet/Neumann boundary conditions outside }  \bf 0:\\
\partial_xu_i(t,R,l)=0,~~ (t,l)\in (0,T)\times(0,K),\\
\forall i\in[\![1,I]\!],~~ u_i(t,x,K)=\psi_i(t,x),~~ (t,x)\in [0,T]\times[0,R],\\
\textbf{Initial condition:}\\
\forall i\in[\![1,I]\!],~~ u_i(0,x,l)=g_i(x,l),~~(x,l)\in[0,R]\times[0,K],\\
\textbf{Continuity condition at } \bf 0:\\
\forall (i,j)\in[\![1,I]\!]^2,~~u_i(t,0,l)=u_j(t,0,l)=u(t,0,l),~~(t,l)\in|0,T]\times[0,K],
\end{cases}
\end{eqnarray}
which has been studied recently by Martinez-Ohavi in \cite{linear PDE}. Before describing the origins and motivations of studying system \eqref{eq PDE 0} - $\mathcal{W}_{alsh}(\mathbb{S})$, let us describe briefly the principal technical characteristics of studying the well-posedness of system \eqref{eq : pde with l}. From a PDE technical aspect, since the variable $l$ drives dynamically the system only at the junction point $\bf 0$ with the presence of the derivative $\partial_lu(t,0,l)$ in the local-time Kirchhoff's boundary condition, the main challenge was to understand the regularity of the solution. Under mild assumptions, it is shown  that classical solutions of the system \eqref{eq : pde with l} belong to the class $\mathcal{C}^{1,2}$ in the interior of each edge and $\mathcal{C}^{0,1}$ in the whole domain (with respect to the time-space variables $(t,x)$). One can expect a regularity in the class $\mathcal{C}^1$ for $l\mapsto u(t,0,l)$ and this is indeed the case (see Theorem 2.4 and point $iv)$ in Definition 2.1 in \cite{linear PDE}). Another technical aspect, was to obtain an H\"{o}lder continuity of the partial functions $l\mapsto \Big(\partial_tu_i(t,x,l),\partial_xu_i(t,x,l),\partial_x^2u_i(t,x,l)\Big)$ for any $x>0$. It is shown using Schauder's estimates, that such regularity is guaranteed by the central assumption on the ellipticity of the diffusion coefficients on each ray, together with the mild dependency of the coefficients and free term with respect to the variable $l$. The comparison theorem holds true in the linear case under such regularity class; see Theorem 2.6 in \cite{linear PDE}.
The reader has to keep in mind that this class of regularity will be naturally used for test functions of system \eqref{eq PDE 0} - $\mathcal{W}_{alsh}(\mathbb{S})$ - since if a viscosity solution is in this regularity class, then the formulation with test functions in this class remains compatible. 
\subsection{Origins and motivations:}
The heat equation posed on networks and its different variants are naturally related to diffusions on graphs introduced in the seminal works of Freidlin and Wentzell \cite{Freidlin} and Freidlin and Sheu \cite{freidlinS}.  

More precisely, given $I$ pairs $(\sigma_{i},b_{i})_{i\in I}$ of mild coefficients of diffusion from $[0,+\infty)$ to $\R$ satisfying the following condition of ellipticity: $\forall i\in [\![1,I]\!],~\sigma_i>0$, and given $\big(\alpha_1,\ldots,\alpha_I)$ positive constants satisfying $\displaystyle \sum_{i=1}^I \alpha_i=1$, it is proved in \cite{Freidlin} that there exists a continuous Feller Markov process $\big(x(\cdot),i(\cdot)\big)$ valued in the star-shaped network, whose generator is given by the following operator:
$$\mathcal{L}:\begin{cases}\mathcal{C}^2(\mathcal{N}_R)\to \mathcal{C}(\mathcal{N}_R),\\
f=f_i(x)\mapsto 
b_i(x)\partial_xf_i(x)+\displaystyle \frac{\sigma_i^2(x)}{2}\partial_x^2f_i(x)\end{cases},$$
with domain
$$D(\mathcal{L}):=\Big\{ f \in \mathcal{C}^2(\mathcal{N}_R),~~\displaystyle \sum_{i=1}^I\alpha_i\partial_xf_i(0)=0\Big\}.$$
Recall that it is shown in \cite{freidlinS} that there exists a one dimensional Wiener process $W$ defined on a probability space $(\Omega,\mathcal{F},\P)$ and adapted to the natural filtration of $\big(x(\cdot),i(\cdot)\big)$, such that the process $\big(x(\cdot)\big)$ satisfies the following stochastic differential equality:
\begin{equation*}
dx(t)= b_{i(t)}(x(t))dt + \sigma_{i(t)}(x(t))dW(t)+d\ell(t) \;,~~0\leq t\leq T.
\end{equation*}
In the above equality, the process $\ell(\cdot)$ is the local time of the process $\big(x(\cdot)\big)$ at the vertex $\bf 0$. Moreover, the process $\ell(\cdot)$ has continuous increasing paths, starts from $0$ and satisfies:
\begin{eqnarray}\label{eq: croissance local time}
\forall t\in [0,T],~~\int_{0}^t\mathbf{1}_{\{x(s)>0\}}d\ell(s)=0,~~\P-\text{a.s.}
\end{eqnarray}
Recall also that the following It\^{o}'s formula was proved in \cite{freidlinS}: 
\begin{eqnarray}\label{Ito Sheu}
\nonumber&\displaystyle df_{i(t)}(x(t))~~=~~ \Big(b_{i(t)}(x(t))\partial_xf_{i(t)}(x(t))+ \frac{1}{2}\sigma_{i(t)}^2(x(t))\partial_{x}^2f_{i(t)}(x(t))\Big)dt+\\
&
\displaystyle \partial_xf_{i(t)}(x(t))\sigma_{i(t)}(x(t))dW(t) 
+ \sum_{i=1}^{I}\alpha_{i}\partial_xf_i(0)d\ell(t),~~\P-\text{a.s,}
\end{eqnarray}
for any sufficiently regular $f$.

From the perspective of stochastic control theory, the question of how to formulate  a problem of stochastic control for this type of processes that satisfies the Itô's rule given in \eqref{Ito Sheu}, is relevant and challenging. As the spider diffusion behaves in each ray like a classical diffusion, the novelty consists exactly in the study of the optimal scattering (or diffraction) of the spider at the junction point $\bf 0$, that is governed by it excursions and the
appropriate and highly non-trivial “Kirchhoff's Law”, appearing in front of the variations of the local time in \eqref{Ito Sheu}:
$$\Big(\sum_{i=1}^{I}\alpha_{i}\partial_xf_i(0)\Big)d\ell(t).$$ 
In this context, let us refer to the recent works in \cite{Kara control}, where an optimal stochastic stopping control problem for a Walsh's planar semi martingale has been studied. Therein, the authors had to face with this technical point, since it is assumed that the process is “immediately dispatched along some ray” when it reaches the vertex point $\bf 0$ (see Section 5). 

Thus, if we have to deal for example with a problem of control with a finite horizon time, it appears clear that the coefficients - especially the diffraction terms $\alpha_i$ - have to depend (first) on the time. Let us now try to collect as much information as possible on the spider at $\bf 0$, to better understand how an eventual possible optimal diffraction behaves. Recall that the local time of the spider is in some way related to the excursions of the spider at $\bf 0$. Indeed, if $\mathcal{N}^\varepsilon(\cdot)$ refers to the number of excursions of size $\varepsilon$ of the spider at $\bf 0$, it is well known from Paul Levy's works called the {\it "Down crossing representation of the local time"}, (see Theorem 2.23, Chapter VI in \cite{Karatzas Book} in the case of a reflected Brownian motion) that we have: 
\begin{align*}
\lim_{\varepsilon\searrow 0}\E\big[~|\varepsilon\mathcal{N}^\varepsilon(t)-\ell(t)|^2~\big]=0.
\end{align*}
Therefore we can suggest making the diffraction coefficients $\alpha_i$ depending on the local time, but that is not all. Let us describe other relevant arguments for this potential dependency. Assume that we face with a controlled spider diffusion, with controls on each ray, namely following an It\^{o}'s rule of type:
\begin{eqnarray}\label{Ito rule control}
\nonumber&\displaystyle df_{i(t)}(x(t))=\Big(b_{i(t)}(x(t),\beta_{i(t)}(t))\partial_xf_{i(t)}(x(t))+ \frac{1}{2}\sigma_{i(t)}^2(x(t),\beta_{i(t)}(t))\partial_{x}^2f_{i(t)}(x(t))\Big)dt\\
&
+\displaystyle \partial_xf_{i(t)}(x(t))\sigma_{i(t)}(x(t),\beta_{i(t)}(t))dW(t) 
+ \sum_{i=1}^{I}\alpha_{i}\partial_xf_i(0)d\ell(t),~~\P-\text{a.s,}
\end{eqnarray}
where the controls $\beta_i(\cdot)$ on each ray are valued on a compact set $\mathcal{B}_i$ of $\R$. It is also legitimate to ask the question of how these controls will intervene in the optimal diffraction of the spider in the neighborhood of the vertex $\bf 0$. For a spider process following \eqref{Ito Sheu}, the local time $\ell$ is related to the second order terms and the coefficients of diffraction $\alpha_i$ by the following approximation (see Remark 2.5 in \cite{freidlinS}) :  
\begin{eqnarray}\label{eq : quadra temps local 1}
\lim_{\varepsilon \searrow 0}~~\mathbb{E}^{\mathbb{P}} \Big[~~\Big|~~\Big(\frac{1}{2\varepsilon}\int_{0}^{t}\mathbf{1}_{\{0\leq x(s)\leq\varepsilon\}}ds\Big)~~-~~\big(\ds \sum_{i=1}^I\frac{\alpha_i}{\sigma^2_i(0)}\big)\ell(t)~~\Big|~~\Big] ~~ =~~0.
\end{eqnarray}
The dependence between the second order plus the diffraction terms and the local time $\ell$ appears then with \eqref{eq : quadra temps local 1}. Recall thus, not only do we notice that the local time increases only when the spider reaches $\bf 0$ via \eqref{eq: croissance local time}, but also that \eqref{eq : quadra temps local 1} implies also that the local time disappears as soon as the terms $\sigma_i(0)$ vanish at $\bf 0$. 
Assume now that the spider follows the controlled dynamic given in \eqref{Ito rule control}, and there exists an optimal control process $\beta_i^*(s)$; depending in the most of the times on the Laplacian and the gradient of a value function $u$ at the interior of each ray, idem:
$$\beta_i^*(s)=\beta_i^*\Big(s,x(s),u\big(x(s)\big),\nabla u\big(x(s)\big),\Delta u\big(x(s)\big)\Big).$$
We can conjecture via the approximation \eqref{eq : quadra temps local 1}, that we will obtain an 'optimal' local time process $l^*$ depending indeed on the controls $\beta_i^*(s)$ in the neighborhood of $\bf 0$. Assume then that coefficients of diffraction $\alpha_i$ depend on the local time and are controlled by a process $\vartheta_s$ valued in some compact set $\mathcal{O}$ of $\R^I$, which should imply that the Kirchhoff's term in \eqref{Ito rule control} is replaced by:
\begin{eqnarray}\label{eq coef controler}
 \sum_{i=1}^I\alpha_i(l(s),\vartheta_s)\partial_xf_i(0).   
\end{eqnarray}
It follows that the optimal stochastic scattering (diffraction) process $\vartheta_s^*$ will depend then on $l^*$ and therefore with a significant chance via \eqref{eq : quadra temps local 1}, also on the behavior of the optimal control terms on each ray: $\beta_i^*(s)$ near $\bf 0$. It appears then natural to consider the local time as an intrinsic variable of the stochastic scattering control problem for the spider, to better understand how a possible optimal diffraction will behave. We will see that this probabilistic intuition - that is indeed a physical interpretation on the behavior of the diffraction of Brownian particles - will lead to certain advancements in term of PDE analysis for Neumann problems with discontinuities.

For this purpose, all the following issues should surely have to be considered:\\
-a) build a spider process having a spinning measure depending on the local time,\\
-b) obtain a comparison theorem for viscosity solution of HJB systems posed on networks, having a non linear Kirchhoff's condition at $\bf 0$,\\
-c) formulate and solve a problem of stochastic optimal scattering control for a spider diffusion, having an optimal spinning measure selected from its own local time.\\
The second problem b) is the main target of this work.\\
Problem a) was achieved recently by Martinez-Ohavi, with the aid of two long and correlated contributions in stochastic analysis and PDE, that are respectively: \cite{Spider} and \cite{linear PDE}.
In \cite{Spider}, the authors take the results stated in \cite{freidlinS} as a starting point and construct 'by hand' a solution of a martingale problem that is purposely designed in order to take the presence of the local time in all the leading coefficients into account. Of course, since the local time is added in the picture, the canonical space has to be extended accordingly. Indeed the martingale problem is formulated as following:\\
-given a time $t\in [0,T)$ (where $T$ refers to a finite horizon time) and a starting point point $\big((x,i),l\big)$ in the 'star-shaped/local-time' state space, the main target is to show that there exists a unique probability measure denoted by $\P_t^{x,i,l}$, defined on the canonical space of continuous maps living on the star-shaped network - times the set of the non negative and non decreasing function - such that for $f$ any regular enough:
$$\big(\mathcal{S}_{pi}-\mathcal{M}_{ar}\big)~~-~~\text{label for the spider martingale problem}:$$
\begin{align} \label{eq:def-V-intro}
&\nonumber\Bigg(f_{i(s)}(s,x(s),l(s))- f_{i}(t,x,l)-\\
&\nonumber\int_{t}^{s}\Big(\partial_tf_{i(u)}(u,x(u),l(u))+\displaystyle\frac{1}{2}\sigma_{i(u)}^2(u,x(u),l(u))\partial_{xx}^2f_{i(u)}(u,x(u),l(u))\\
&\nonumber\displaystyle+b_{i(u)}(u,x(u),l(u))\partial_xf_{i(u)}(u,x(u),l(u))\Big)du-\\
&\int_{t}^{s}\big(\partial_{l}f(u,0,l(u))+\displaystyle\sum_{j=1}^{I}\alpha_j(u,l(u))\partial_{x}f_{j}(u,0,l(u))\big)dl(u)~\Bigg)_{t\leq s\leq T},
\end{align}
is a martingale under the probability measure $\P_t^{x,i,l}$ and the natural filtration generated by the canonical process $(x(s),i(s),l(s))_{s\in [t,T]}$. Here $(l(s))_{s\in [t,T]}$ stands for the local time of the expected spider process at the junction point $\bf 0$. Because $l(\cdot)$ has continuous increasing paths, and increases only when the spider reaches $\bf 0$, it appears that the following {\it local-time Kirchhoff's transmission condition}:
\begin{equation}
\label{eq:condition-transmission-intro}
\partial_l f(t,0,l)+\displaystyle \sum_{i=1}^I \alpha_i(t,l)\partial_x f_i(t,0,l)=0,\;\;\;(t,l)\in(0,T)\times(0,+\infty)
\end{equation}
must hold in the martingale formulation in the front of the variations of $l$, for any $f$ regular enough. 
The main - Theorem 3.1 in \cite{Spider} - states that the martingale problem $\big(\mathcal{S}_{pi}-\mathcal{M}_{ar}\big)$ is well-posed on the corresponding canonical space.
The existence proof relies on a very careful adaptation of the seminal construction for solutions of classical martingale problems that have $\R^d$ as the underlying state space, and combines concatenation of probability measures with tension arguments. For more details on the choice of this construction, we refer the reader to the Introduction in \cite{Spider}. At last, let us quote that \cite{Spider} provides the first result in literature of existence of a Walsh processes having non constant spinning measure.

As pointed before, the authors had to face with the problem of studying the well-posedness of the corresponding parabolic operator given in \eqref{eq : pde with l}, and this was achieved in \cite{linear PDE}. The results contained in \cite{linear PDE} are of crucial importance when turning to the difficult problem of uniqueness for the martingale problem $\big(\mathcal{S}_{pi}-\mathcal{M}_{ar}\big)$. Note that the results contained in \cite{linear PDE} extend and also improve (see Section 3 in \cite{linear PDE}) those obtained by Von Below in \cite{Von Below}, which were -- up to our knowledge -- the only reference on the well-posedness of heat equations on graphs with time variable coefficients and classical Kirchhoff's condition:
\begin{eqnarray*}
\sum_{i=1}^I \alpha_i(t)\partial_x f_i(t,0)=0.
\end{eqnarray*}
 
Finally, to conclude this Sub section, let us make the link with control theory. Assume that we face with a stochastic control problem formulated in the weak sense for probability measures solving the martingale problem $\big(\mathcal{S}_{pi}-\mathcal{M}_{ar}\big)$ given in \eqref{eq:def-V-intro}. We introduce for this purpose controlled processes $(\beta_{i(u)})_{u\in [t,T]}$ on each ray $\mathcal{R}_i$ and $(\vartheta_u)_{u\in [t,T]}$ at the vertex $\bf 0$, now progressively measurable with respect to the filtration of the canonical process $(x(u),i(u),l(u))_{u\in [t,T]}$. Consider the set of all probability measures - denoted by $\mathcal{A}(t,x,i,l)$ -  such that for any test function $f$:
\begin{align} \label{Martingale control}
&\nonumber\Bigg(f_{i(s)}(s,x(s),l(s))- f_{i}(t,x,l)-\\
&\nonumber\int_{t}^{s}\Big(\partial_tf_{i(u)}(u,x(u),l(u))+\displaystyle\frac{1}{2}\sigma_{i(u)}^2(u,x(u),l(u),\beta_{i(u)})\partial_{xx}^2f_{i(u)}(u,x(u),l(u))\\
&\nonumber\displaystyle+b_{i(u)}(u,x(u),l(u),\beta_{i(u)})\partial_xf_{i(u)}(u,x(u),l(u))\Big)du-\\
&\int_{t}^{s}\big(\partial_{l}f(u,0,l(u))+\displaystyle\sum_{j=1}^{I}\alpha_j(u,l(u),\vartheta_u)\partial_{x}f_{j}(u,0,l(u))\big)dl(u)~\Bigg)_{t\leq s\leq T},
\end{align}
is a martingale under the natural assumptions. Assume that an agent deals with a problem of optimization, having a cost $h_0$ at $\bf 0$ and $h_i$ on each ray $\mathcal{R}_i$. Fix now $t=0$, and assume that all the coefficients do not depend on the current time anymore. Define the following cost function: 
\begin{eqnarray}\label{eq : cost function}
w:=\begin{cases}
\mathcal{A}(t=0,x,i,l) \to \R \\
\P_0^{x,i,l}\mapsto\displaystyle\E^{\P_0^{x,i,l}}\Big[~~\int_{0}^{+\infty} e^{-\lambda u}h_{i(u)}(x(u),l(u),\beta_{i(u)})du~~+\\
\hspace{2.5cm}\displaystyle \int_{0}^{+\infty}e^{-\lambda u}h_0(l(u),\vartheta_u)dl(u)~~\Big]
\end{cases},
\end{eqnarray}
and for instance a corresponding value function $u$ by:
\begin{eqnarray}\label{eq : value function}
u:=\begin{cases}
\mathcal{N}\times[0,+\infty)\to \R \\
\big(x,i,l\big)\mapsto \sup \Big\{~~w\big(\P_0^{x,i,l}\big),~~\P^{x,i,l}_0\in\mathcal{A}\big(t=0,x,i,l\big)~~\Big\} 
\end{cases}.
\end{eqnarray}
It is expected that the value function $u$ will be a viscosity solution of the following system:
\begin{eqnarray}\label{HJB for the control}
\begin{cases}
\textbf{HJB equation parameterized by the local-time on each ray:}\\ 
-\lambda u_i(x,l)+\underset{\beta_i\in \mathcal{B}_i}{\sup}\Big\{\ds \frac{1}{2}\sigma_i(x,l,\beta_i)^2\partial^2_xu_i(x,l)+
b_i(x,l,\beta_i)\partial_xu_i(x,l)+\\
h_i(x,l,\beta_i)\Big\}=0,~~(x,l)\in(0,+\infty)^2,\\
\textbf{Non linear local-time Kirchhoff's boundary condition at } \bf 0:\\
\partial_lu(0,l)+\underset{ \vartheta \in \mathcal{O}}{\sup} \Big\{\ds\sum_{i=1}^I\alpha_i(l,\vartheta)\partial_xu_i(0,l)+h_0(l,\vartheta)\Big\}=0,~~l\in(0,+\infty),\\
\textbf{Continuity condition at } \bf 0:\\ 
\forall (i,j)\in[\![1,I]\!]^2,~~\forall l\in [0,+\infty),~~u_i(0,l)=u_j(0,l).
\end{cases}
\end{eqnarray}
We thus obtain the origin of the study of a system of type \eqref{eq PDE 0} - $\mathcal{W}_{alsh}(\mathbb{S})$. Using some basic modifications, we see that its existence is a consequence of the control theory. Note that in order to obtain the corresponding Dirichlet condition in \eqref{eq PDE 0}, it is enough to introduce the corresponding exit times of the processes $x$ and $l$, and modify the cost function appropriately.
Since this contribution is quite long, we have decided to push back in an future work the study of the problem of control. Therein, we will also try to discuss on the possible applications in the field of quantum physics. Indeed, the study of the optimal diffraction for the Walsh diffusion variety, has important implications when one tries to study the diffusive behavior of particles subjected to scattering (or diffraction), for which very little physical, understanding currently exists. The theory of quantum trajectories states that quantum systems can be modelled as scattering processes, and we refer the reader to \cite{Scatterin Theory}, for more appropriate details. As a consequence, problems of optimal light scattering have attracted several scientists, for their importance in advanced photonics technologies, as in on-chip interconnects, bioimaging, solar-cells, heat-assisted magnetic recording, and many orders.
\subsection{On the terminology of the system \eqref{eq PDE 0} and its possible extensions:}  
We start by explaining why we decided to call the system given in \eqref{eq PDE 0}: {\it Walsh's spider Hamilton-Jacobi-Bellman system} - $\mathcal{W}_{alsh}(\mathbb{S})$.

System \eqref{eq PDE 0} - $\mathcal{W}_{alsh}(\mathbb{S})$ and its variants, whether in the linear or non-linear frameworks, are more often called in the literature: ''system of PDE posed on a star-shaped network''. See for instance amongst the vast works on this subject:  \cite{control 5}, \cite{Lions Souganidis 1}, \cite{Lions Souganidis 2}, \cite{Lions Souganidis 3}, \cite{linear PDE}, \cite{Ohavi PDE}, \cite{Von Below}, \cite{Below 3}, \cite{Below 4}...

It is crucial to note that in the most of the works in literature dealing with the analysis of these PDE systems posed on networks, the angles appearing in the geometry of the network do not play any role in obtaining existence, regularity, comparison theorem... Indeed it appears that we can consider any star-shaped network with exactly the same rays in any Euclidean space $\R^n,~n\ge 2$; spaced by the same number of angles, without modifying the mathematical analysis of these PDE systems.
Indeed diffusions on graphs can also be seen as part of the family of Walsh processes introduced in the epilogue \cite{Walsh}, let us describe briefly this point. Assume now that we are given $I$ distinct angles $(\theta_1,\ldots,\theta_{I})\in [0,2\pi)^I$ that parameterize uniquely the star shaped network (belonging to the plan $\R^2)$:
$$\displaystyle \mathcal{N}_R=\bigcup_{i=1}^I\mathcal{R}_i,$$
in the following sense:
$$\forall i\in [\![1,I]\!],~~(\vec{\bf{0} \mathcal{R}_i},\vec{\bf{0} \mathcal{R}_{i+1}})=\theta_i,$$
(with the convention $\mathcal{R}_{I+1}=\mathcal{R}_1$).
Consider now a 'spider' probability measure $\mathbb{S}$ on $[0,2\pi)$ with constant support $(\theta_1,\ldots,\theta_{I})$, such that:
$$\forall i\in[\![1,I]\!],~~\mathbb{S}(\theta_i)=\alpha_i,~~\sum_{i=1}^I\alpha_i=1.$$
We see then that the Kirchhoff's term appearing in the It\^{o}'s formula \eqref{Ito Sheu}, can be rewritten as in the classical formulation used for Walsh processes in the literature by:
$$\sum_{i=1}^{I}\alpha_{i}\partial_xf_i(0)d\ell(t)=\Big(\int_0^{2\pi}\partial_xf_{\theta}(0)\mathbb{S}(d\theta)\Big)d\ell(t),$$
(where abusively we have denoted $\partial_xf_i(0)=\partial_xf_{\theta_i}(0)$).
We now see  under this formulation that the geometry of the star-shaped network, and consequently in some way a PDE system posed on this structure, is uniquely characterized by the support of the spider measure $\mathbb{S}$, namely the angles $(\theta_1,...,\theta_I)$. This is the main reason which pushed us to call the system \eqref{eq PDE 0}: Walsh's spider Hamilton-Jacobi-Bellman system - $\mathcal{W}_{alsh}(\mathbb{S})$. We have in mind consequently the extension of our results, from the special case when $\mathbb{S}$ has a constant support, to more general angle measures.

\subsection{Main contributions and novelties:}
First of all, it should be noted that the supremums and infinimums appearing in the system  \eqref{eq PDE 0} do not play a preponderant role in the proof of our comparison theorem, and no convexity is required. By adaptation of the main ideas of this work, the reader can easily obtain the same corresponding results for other situations, for instance other kinds of supremums or infimums, such as the Isaacs HJB equation. We wanted to formulate the non linearity using Hamiltonians that appear naturally in stochastic control theory. The intuition is for the reader to keep an eye on the underlying process, in order to better understand the behavior of the process at the junction point and its non-stickiness. It is more precisely the understanding of the non-stickiness, in particular its proof in the case of a Walsh's spider diffusion, that gave us some intuitions for the construction of the test functions.

Fully non linear Hamiltonians can be considered, as soon as they satisfy an ellipticity condition on each domain $\mathcal{R}_i=[0,R]\times \{i\},~i\in[\![1,I]\!]$; quadratic growth with respect to the gradient, whereas the Kirchhoff's condition at the vertex $\bf 0$ 
must simply be strictly increasing with respect to the gradient. We aim in an upcoming work, to obtain existence and uniqueness of system \eqref{eq PDE 0} with more general Hamiltonians, and especially as explained in the last Sub section, with a more general angular measure. Let us finally mention that our results can be also extended to time-dependent problems using classical arguments arising from the theory of viscosity solutions.

Another key point when we study the behavior of a Walsh's spider motion, is to notice that all discontinuities generated by the coefficients $\big(b_i(0)\neq b_j(0),\sigma_i(0)\neq \sigma_j(0)\big)_{i\neq j}$, disappear in the It\^{o}'s rule \eqref{Ito Sheu} with the aid of what we call the {\it non-stickiness} condition satisfied by the Walsh's spider $\big(x(t),i(t)\big)_{t\ge 0}$ at $\bf 0$, that reads:
$$\forall t \ge 0,~~\int_0^t\mathbf{1}_{\{x(s)=0\}}ds=0,~~\P~~a.s.$$
To obtain this {\it non-stickiness} condition, the authors in \cite{freidlinS} used an ordinary differential equation (ODE). This ODE involves coefficients depending on the speed measure $s_i(x)dx$ of the Walsh's spider on each ray, defined by:
\begin{eqnarray}\label{eq speed 1}
  \forall x>0,~~\forall i \in[\![1,I]\!],~~s_i(x)dx=\frac{2}{\sigma^2_i(x)}\exp\big(\int_0^x\frac{2b_i(z)}{\sigma^2_i(z)}dz\big)dx.  
\end{eqnarray}
This fact has given us the intuitions to build test functions of system \eqref{eq PDE 0} at the vertex $\bf 0$.
More precisely, if one considers the following Hamiltonians $\big(H_i,~i\in[\![1,I]\!]\big)$:
\begin{align*}
H_i:=\begin{cases}
    [0,R]\times[0,K]\times \R^3\to \R,\\
    (x,l,u,p,S)\mapsto \lambda u+\underset{\beta_i \in \mathcal{B}_i}{\sup}\Big\{\ds-\sigma_i(x,l,\beta_i)S+b_i(x,l,\beta_i)p+h_i(x,l,\beta_i)\Big\},
\end{cases}   
\end{align*}
appearing on each ray of system $\eqref{eq PDE 0}$, we can define similarly to the speed measure $s_i(x)dx$, the  
{\it speed of the Hamiltonians}  $\mathcal{S}_{peed}\big((H_i)_{i\in[\![1,I]\!]}\big)$ as: 
\begin{align*}
\mathcal{S}_{peed}\big((H_i)_{i\in[\![1,I]\!]}\big):=\begin{cases}
    [0,R]\times[0,K]\times \R\to \R,\\
    (x,l,p)\mapsto \underset{\beta_i \in \mathcal{B}_i}{\sup}\Big\{\ds \frac{b_i(x,l,\beta_i)p+h_i(x,l,\beta_i)}{\sigma_i(x,l,\beta_i)}\Big\}
\end{cases}.
\end{align*}
To prove the comparison theorem for system \eqref{eq PDE 0}, we will build test functions at the junction point $\bf 0$ solutions of ODE with coefficients that may be viewed as a kind of envelope of all possible errors of the {\it speed of the Hamiltonians}  
$\mathcal{S}_{peed}\big((H_i)_{i\in[\![1,I]\!]}\big)$. The key point in this construction is to impose a local time derivative with respect to variable $l$, at $\bf 0$ - $\partial_l\phi(0,l)$ - that will absorb all the error term induced by  - that we have decide to call - the {\it Kirchhoff's speed of the Hamiltonians} - $\mathcal{K}_{f}\mathcal{S}_{peed}\big((H_i)_{i\in[\![1,I]\!]}\big)$ - defined by:
\begin{align*}
\mathcal{K}_{f}\mathcal{S}_{peed}\big((H_i)_{i\in[\![1,I]\!]}\big):=\begin{cases}
    [0,R]\times[0,K]\times \R\to \R,\\
    (x,l,p)\mapsto \ds \inf_{\vartheta\in \mathcal{O}}\Big\{\sum_{i=1}^I\mathbb{S}_i(l,\vartheta)\mathcal{S}_{peed}(H_i)(x,l,p)\Big\}
\end{cases}.
\end{align*}

Up to our knowledge, this is the first result of uniqueness for HJB elliptic PDE system posed on a star shaped-network, having a non linear Kirchhoff's condition and non vanishing viscosity at the vertex.

Recall that in the theory of viscosity solutions, the formulation of upper and sub solutions for a second order non linear problem with Neumann $N$ (or Kirchhoff"s) boundary condition:
\begin{eqnarray}\label{eq Neumann}
\begin{cases}
    H\big(x,u(x),\nabla u(x),\Delta u(x)\big)=0,~~x\in \Omega,\\
    N\big(x,u(x),\nabla u(x)\big)=0,~~x\in \partial \Omega
\end{cases},    
\end{eqnarray}
(where $\Omega$ denotes a smooth open set of $\R^n$) can be stated in the strong sense namely using only the Neumann condition $N(\cdot)$ at the boundary for both upper and sub solutions, or in the weak sense with both terms $N(\cdot)$ and $H(\cdot)$, considering at the boundary $\partial \Omega$ the term $N(\cdot)\wedge H(\cdot)$ for sub solutions and $N(\cdot)\vee H(\cdot)$ for super solutions (see for instance Section 7 in \cite{User guide}). There are already a lot of comparison and existence results for viscosity solutions of second order PDEs with general Neumann type boundary conditions. We refer for this to \cite{Barles Neuman 1}, \cite{Barles Neuman 2}, \cite{Daniel game},   \cite{Hishi neumann}, \cite{Lions neumann}, \cite{Sato} and references therein.
Recall that in almost all the cases the proof of the comparison theorem consists in introducing the doubling variable method, with the function:
$$\forall \varepsilon>0,~~w_\varepsilon\big(x,y)=u(x)-v(y)-\frac{|x-y|_{\R^n}^2}{2\varepsilon^2},~~(x,y)\in \R^n,$$
(where $u$ is a super solution, whereas $v$ is a sub solution of \eqref{eq Neumann}) and to obtain a contradiction, passing to the limit when $\varepsilon \searrow 0$, locally at any point $x\in \overline{\Omega}$. When $x\in \partial \Omega$, both conditions $F$ and $H$ are then considered, and the Hamiltonian $H$ has to be continuous at the neighborhood of any point $x\in \partial \Omega$. Consequently, we will get that this method will fail for system \eqref{eq PDE 0}, because of the discontinuities of the Hamiltonians at $\bf 0$.

We will see that the construction of the test function for proving the comparison theorem in this work, allow us to formulate super and sub viscosity solution in the strong sense at $\bf 0$. Note that the {\it non linear local-time Kirchhoff's boundary transmission} appears at $\bf 0$, without considering any values of the Hamiltonians at the vertex. This allows us to avoid all the discontinuities induced by the Hamiltonians at $\bf 0$ and to get over this technical point.

Finally, as a consequence of the main results of this contribution, we are able to obtain a comparison theorem in the case when the dependency with respect to the variable $l$ is removed. More precisely, using this new technique consisting of adding a variable representing the local time at $\bf 0$, we are able to obtain a comparison theorem for HJB systems posed on star-shaped networks, with non linear Kirchhoff's boundary condition at $\bf 0$, and no vanishing viscosity at the vertex. Moreover, this comparison theorem is stated with a Neumann (here Kirchhoff's) boundary condition in the strong sense, namely without any dependency of the Hamiltonians at $\bf 0$, that is - up to our knowledge - a new advance in the theory of non linear PDE problems having Neumann's boundary conditions. We hope that the new techniques introduced in our contribution can improve the analysis of Neumann's problems posed on $\R^n$, especially the case of PDE's systems connected on manifolds with discontinuous coefficients.
\subsection{Review of literature:}
To finish this Introduction, let us give an account of the main works that have been done in literature for similar systems close to \eqref{eq PDE 0}. Note that the key fact in the most all of them, is to consider a vanishing viscosity at the vertex $\bf 0$, which is not the case here. 

We first refer to the recent monograph \cite{Barles book}, that presents the most recent developments in the study of Hamilton-Jacobi Equations and control problems with discontinuities (see Part III for the case of problems on networks).
In \cite{Lions Souganidis 1}, the authors introduce a notion of state-constraint viscosity solutions for one dimensional “junction” - type problems for first order Hamilton-Jacobi equations with non-convex coercive Hamiltonians and study its well-posedness and stability properties. Let us quote that in this work, the main results do not require any convexity conditions on the Hamiltonians, contrary to all the previous literature that is based on the control (deterministic) theoretical interpretation of the problem. Among the long list of references on this topic with convex Hamiltonians, we can cite for instance: \cite{control 1}, \cite{control 2}, \cite{control 3}, \cite{control 4}, \cite{control 5}, \cite{control 6}. For recent works on systems of conservative laws posed on junctions, we refer also to \cite{Carda junction 1} and \cite{Carda junction 2} with the references therein. 
In \cite{Lions Souganidis 2}, the authors have studied multi-dimensional junction problems for first and second-order PDE with Kirchhoff-type
Neumann boundary conditions, showing that their generalized viscosity solutions are unique, but still with a vanishing viscosity at the vertex for the second order terms. Finally, let us cite the interesting approach studied in \cite{Lions Souganidis 3}, where it is considered star-shaped tubular domains consisting of a number of non-intersecting, semi-infinite
strips of small thickness that are connected by a central region.  It is shown that classical regular solutions of uniformly elliptic partial differential equations converge in the thin-domain limit, to the unique solution of a second-order partial differential equation on the network satisfying an effective Kirchhoff-type transmission condition at the junction.
\subsection{Organisation of the paper:}
The paper is organized as follows: In Section \ref{sec : Notations}, we introduce the main notations and we state the main Theorem \ref{th: compa theorem principale} of this work. Note that Section \ref{sec : Notations} contains also a comparison theorem, Theorem \ref{theo compa normal} of system of type \eqref{eq PDE 0} when the local-time variable $l$ is removed, inducing then the 'classical' non linear Kirchhoff's boundary condition already investigated in the literature. Section \ref{sec : exemple} is dedicated to give a sketch of proof in the linear case of our method, for simple elliptic PDE, in order to give to the reader some intuitions and ideas that will be used for the construction of tests functions for our central Theorem \ref{th: compa theorem principale}. Finally in Section \ref{sec preuve theorem principal}, we prove our main Theorem \ref{th: compa theorem principale}, that is the comparison principle for system \eqref{eq PDE 0}.

\section{Notations and Definitions}\label{sec : Notations}
In this Section, we introduce the main notations-definitions and we state our main Theorem \ref{th: compa theorem principale}.
In all this work, we fix $R>0$ and $K>0$ the boundaries of the 'space/local-time' domain
$$(0,R)\times (0,K)\ni(x,l)$$
where the system \eqref{eq PDE 0} will be studied.

Let $I\in \mathbb{N}^*$ be the number of edges and $R>0$ be the common length of each ray. The bounded star-shaped compact network $\mathcal{N}_R$ is defined by:
$$\mathcal{N}_R=\bigcup_{i =1}^I \mathcal{R}_{i}$$
where
\begin{eqnarray*} \forall i\in [\![1,I]\!]~~\mathcal{R}_{i}:=[0,R]~~\text{and}~~\forall (i,j)\in [\![1,I]\!]^2,~~i\neq j,~~\mathcal{R}_{i}\cap \mathcal{R}_{j}=\bf 0.
\end{eqnarray*}
The intersection of all the rays $(\mathcal{R}_{i})_{1 \leq i\leq I}$ is called the junction point and is denoted by $\bf 0$.\\
We identify  all the points of $\mathcal{N}_R$ by couples $(x,i)$ (with $i \in[\![1,I]\!], x\in|0,R]$), such that we have: $(x,i)\in \mathcal{N}_R$, if and only if $x\in \mathcal{R}_{i}$.

Let $(\mathcal{B}_i)_{i\in[\![1,I]\!]}$ be a collection of $I$ compact sets of $\R$ and $\mathcal{O}$ a compact set of $\R^I$. We introduce the following data:
$$\textbf{Data:~} (\mathcal{D})~~\begin{cases}
\Big(\sigma_i \in \mathcal{C}\big([0,R]\times [0,K]\times \mathcal{B}_i,\R\big)\Big)_{i\in[\![1,I]\!]}\\
\Big(b_i \in \mathcal{C}\big([0,R]\times[0,K]\times\mathcal{B}_i,\R\big)\Big)_{i\in[\![1,I]\!]}\\
\Big(h_i \in \mathcal{C}\big([0,R]\times [0,K]\times \mathcal{B}_i,\R\big)\Big)_{i\in[\![1,I]\!]}\\
\Big(\mathbb{S}_i\in \mathcal{C}\big([0,K]\times \mathcal{O},\R\big)\Big)_{i\in[\![1,I]\!]}\\
h_0\in \mathcal{C}\big([0,K]\times \mathcal{O},\R\big),\\
\Big(\chi_i \in \mathcal{C}\big([0,K],\R\big)\Big)_{i\in[\![1,I]\!]},\\
\Big(\mathcal{T}_i \in \mathcal{C}\big([0,R],\R\big)\Big)_{i\in[\![1,I]\!]}.
\end{cases}.
$$
We assume that the data $(\mathcal{D})$ satisfy the following assumptions: (where $({\bf S})$ stands for the coefficients of diffraction for the spider, $({\bf E})$ for ellipticity, and $({\bf R})$ for Lipschitz regularity uniformly with respect to the control variables):
\newpage
$$\textbf{Assumption } (\mathcal{H})$$
\begin{align*}
&({\bf S})~~\exists\;\ds \underline{\zeta}>0,~~\forall i \in [\![1,I]\!], ~~\forall (l,\vartheta)\in [0,K]\times\mathcal{O},~~\mathbb{S}_i(l,\vartheta)~~\ge~~\underline{\zeta}.\\
&({\bf E})~~\exists\,\,\underline{\sigma} >0,~~\forall i \in [\![1,I]\!], ~~\forall (x,l,\beta_i)\in [0,R]\times[0,K]\times \mathcal{B}_i,~~\sigma_i(x,l,\beta_i)~~\ge~~\underline{\sigma}.\\
&({\bf R})~~\exists\,\,(|b|,|h|,\overline{\zeta},\overline{\sigma})\in (0,+\infty)^4,~~\forall i \in [\![1,I]\!],\\
&({\bf R} - i)\hspace{0,5 cm}\displaystyle\sup_{x,l,\beta_i}|b_i(x,l,\beta_i)|\;\;\,+\sup_{l,\beta_i}\sup_{(x,y),\;x\neq y} \frac{|b_i(x,l,\beta_i)-b_i(y,l,\beta_i)|}{|x-y|}\\
&\displaystyle\hspace{3,5cm}+\sup_{x,\beta_i}\sup_{(l,l'),\;l\neq l'} \frac{|b_i(x,l,\beta_i)-b_i(x,l',\beta_i)|}{|l-l'|}~~\leq~~|b|,\\
&({\bf R} - ii)\hspace{0,5 cm}\displaystyle\sup_{x,l,\beta_i}|\sigma_i(x,l,\beta_i)|\;\;\,+\sup_{l,\beta_i}\sup_{(x,y),\;x\neq y} \frac{|\sigma_i(x,l,\beta_i)-\sigma_i(y,l,\beta_i)|}{|x-y|}\\
&\displaystyle\hspace{3,5cm}+\sup_{x,\beta_i}\sup_{(l,l'),\;l\neq l'} \frac{|\sigma_i(x,l,\beta_i)-\sigma_i(x,l',\beta_i)|}{|l-l'|}~~\leq~~\overline{\sigma},\\
&({\bf R} - iii)\hspace{0,5 cm}\displaystyle\sup_{x,l,\beta_i}|h_i(x,l,\beta_i)|\;\;\,+\sup_{l,\beta_i}\sup_{(x,y),\;x\neq y} \frac{|h_i(x,l,\beta_i)-h_i(y,l,\beta_i)|}{|x-y|}\\
&\displaystyle\hspace{3,5cm}+\sup_{x,\beta_i}\sup_{(l,l'),\;l\neq l'} \frac{|h_i(x,l,\beta_i)-h_i(x,l',\beta_i)|}{|l-l'|}~~\leq~~|h|,\\
&({\bf R}- iv)\hspace{0,5 cm}\displaystyle \displaystyle\sup_{l,\vartheta}|\mathbb{S}_i(l,\vartheta)|+\sup_\vartheta\sup_{(l,l'),\;l\neq l'} \frac{|\mathbb{S}_i(l,\vartheta)-\mathbb{S}_i(l',\vartheta)|}{|l-l'|}~~\leq ~~\overline{\zeta},\\
&({\bf R}- v)\hspace{0,5 cm}\displaystyle \displaystyle\sup_{l,\vartheta}|h_0(l,\vartheta)|+\sup_\vartheta\sup_{(l,l'),\;l\neq l'} \frac{|h_0(l,\vartheta)-h_0(l',\vartheta)|}{|l-l'|}~~\leq ~~|h|.
\end{align*}
\vspace{0,3 cm}
\textbf{Fix now all around this work: $\lambda>0$.}

In this study, we will obtain a comparison theorem for the following Walsh's spider Hamilton-Jacobi-Bellman system $\mathcal{W}_{alsh}(\mathbb{S})$, having a {\it non linear local-time Kirchhoff's boundary transmission at $\bf 0$}:
\begin{eqnarray}\label{eq system Walsh}
\nonumber &\mathcal{W}_{alsh}(\mathbb{S}):=\\
&\begin{cases}
\textbf{HJB equation parameterized by the local-time on each ray:}\\ 
\lambda u_i(x,l)+\underset{\beta_i\in \mathcal{B}_i}{\sup}\Big\{-\sigma_i(x,l,\beta_i)\partial^2_xu_i(x,l)+
b_i(x,l,\beta_i)\partial_xu_i(x,l)+\\
h_i(x,l,\beta_i)\Big\}=0,~~(x,l)\in(0,R)\times (0,K),\\
\textbf{Non linear local-time Kirchhoff's boundary condition at } \bf 0:\\
\partial_lu(0,l)+\underset{ \vartheta \in \mathcal{O}}{\inf} \Big\{\ds\sum_{i=1}^I\mathbb{S}_i(l,\vartheta)\partial_xu_i(0,l)+h_0(l,\vartheta)\Big\}=0,~~l\in(0,K),\\
\textbf{Dirichlet boundary conditions outside } \bf 0:\\ 
u_i(R,l)=\chi_i(l),~~l\in[0,K],\\
u_i(x,K)=\mathcal{T}_i(x),~~x\in[0,R],\\
\textbf{Continuity condition at } \bf 0:\\ 
\forall (i,j)\in[\![1,I]\!]^2,~~\forall l\in [0,K],~~u_i(0,l)=u_j(0,l).
\end{cases}
\end{eqnarray}
Given the following Hamiltonians $\big(H_i,~i\in[\![1,I]\!]\big)$ defined by:
\begin{align*}
H_i:=\begin{cases}
    [0,R]\times[0,K]\times \R^3\to \R,\\
    (x,l,u,p,S)\mapsto \underset{\beta_i \in \mathcal{B}_i}{\sup}\Big\{\ds\lambda u-\sigma_i(x,l,\beta_i)S+b_i(x,l,\beta_i)p+h_i(x,l,\beta_i)\Big\}.
\end{cases}   
\end{align*}
we will often refer in this work to the 
{\it speed of the Hamiltonians} $\big(\mathcal{S}_{peed}\big((H_i)_{i\in[\![1,I]\!]}\big)\big)$ and the {\it Kirchhoff's speed of the Hamiltonians} $\big(\mathcal{K}_{f}\mathcal{S}_{peed}\big((H_i)_{i\in[\![1,I]\!]}\big)\big)$, that are both defined by:
\begin{align*}
\mathcal{S}_{peed}(H_i):=\begin{cases}
    [0,R]\times[0,K]\times \R\to \R,\\
    (x,l,p)\mapsto \underset{\beta_i \in \mathcal{B}_i}{\sup}\Big\{\ds \frac{b_i(x,l,\beta_i)p+h_i(x,l,\beta_i)}{\sigma_i(x,l,\beta_i)}\Big\}
\end{cases},
\end{align*}
and
\begin{align*}
\mathcal{K}_{f}\mathcal{S}_{peed}\big((H_i)_{i\in[\![1,I]\!]}\big):=\begin{cases}
    [0,R]\times[0,K]\times \R\to \R,\\
    (x,l,p)\mapsto \ds \inf_{\vartheta \in \mathcal{O}}\Big\{\sum_{i=1}^I\mathbb{S}_i(l,\vartheta)\mathcal{S}_{peed}(H_i)(x,l,p)\Big\}
\end{cases}.
\end{align*}
For a given sequence of real numbers $(u_{\varepsilon_1,\ldots,\varepsilon_p})$ indexed by $p$ ($p\in \mathbb{N}^*$) variables $(\varepsilon_1,\ldots,\varepsilon_p)$, we will denote (if the limit exists):
\begin{eqnarray*}
\limsup_{\underset{\varepsilon_j\to a_j,\ldots,\varepsilon_p \to a_p}{\varepsilon_1\to a_1,...,\varepsilon_{j-1}\to a_{j-1},}}u_{\varepsilon_1,\ldots,\varepsilon_p}=\limsup_{\varepsilon_1\to a_1}\ldots\limsup_{\varepsilon_j\to a_j}\ldots\limsup_{\varepsilon_p\to a_p} u_{\varepsilon_1,\ldots,\varepsilon_p},
\end{eqnarray*}
where $j\in [\![2,p]\!]$ and $(a_1,\ldots,a_p)\in \overline{\R}^p$.
In order to remain consistent with the results obtained in \cite{linear PDE}, more precisely with the class of regularity of the solutions in the linear framework (see Introduction in \cite{linear PDE}), we introduce the following space of test functions for continuous viscosity solutions of the system $\mathcal{W}_{alsh}(\mathbb{S})$:
\begin{align*}
&\mathcal{C}^{2,0}_{{\bf 0},1}\big(\mathcal{N}_R\times [0,K]\big):=\Big\{f:\mathcal{N}_R\times [0,K],~~((x,i),l)\mapsto f_i(x,l)~\Big\vert\\
&\hspace{2 cm}~\forall i\in [\![1,I]\!], \;f_i:[0,R]\times [0,K]\to\R,\,(x,l)\mapsto f_i(x,l)\in \mathcal{C}^{2,0}([0,R]\times [0,K]),\\
&\hspace{2 cm}~ \forall (i,j,l)\in [\![1,I]\!]^{2}\times [0,K], \,f_i(0,l)=f_j(0,l)=f(0,l),
\\&\hspace{2 cm}~f(0,\cdot)\in \mathcal{C}^{1}\big([0,K]\big)\Big\},
\end{align*}
where $\mathcal{C}^{2,0}([0,R]\times [0,K])$ states for the set of functions that are $\mathcal{C}^2$ w.r.t. the first variable and continuous w.r.t. the second one.\\
We continue this Section by giving the definition of continuous viscosity super and sub solutions that belong to $\mathcal{C}\big(\mathcal{N}_R\times [0,K]\big)$, defined by:
\begin{align*}
&\mathcal{C}\big(\mathcal{N}_R\times [0,K]\big):=\Big\{f:\mathcal{N}_R\times [0,K],~~((x,i),l)\mapsto f_i(x,l)~\Big\vert\\
&\hspace{2 cm}~\forall i\in [\![1,I]\!], \;f_i:[0,R]\times [0,K]\to\R,\,(x,l)\mapsto f_i(x,l)\in \mathcal{C}^{0}([0,R]\times [0,K]),\\
&\hspace{2 cm}~ \forall (i,j,l)\in [\![1,I]\!]^{2}\times [0,K], \,f_i(0,l)=f_j(0,l)=f(0,l) \Big\}.
\end{align*} for the Walsh's spider HJB system - $\mathcal{W}_{alsh}(\mathbb{S})$ - given in \eqref{eq system Walsh}.
\begin{Definition}\label{def viscosity}
Let $u\in \mathcal{C}\big(\mathcal{N}_R\times [0,K]\big)$.

a) We say that $u$ is a continuous viscosity super solution of the $\mathcal{W}_{alsh}(\mathbb{S})$ system \eqref{eq system Walsh}, if for all test function $\phi \in \mathcal{C}^{2,0}_{{\bf 0},1}\big(\mathcal{N}_R\times [0,K]\big)$ and for all local minimum point 
$(x_\star,i_\star,l_\star)\in [0,R]\times [\![1,I]\!] \times [0,K]$ of $u-\phi$, with $(u-\phi)_{i_\star}(x_\star,l_\star)=0$, we have:
\begin{align*}
\begin{cases}
\lambda \phi_{i_\star}(x_\star,l_\star)+\underset{\beta_{i_\star}\in \mathcal{B}_{i_\star}}{\sup}\Big\{-\sigma_{i_\star}(x_\star,l_\star,\beta_i)\partial^2_x\phi_{i_\star}(x_\star,l_\star)+\\
b_{i_\star}(x_\star,l_\star,\beta_i)\partial_x\phi_{i_\star}(x_\star,l_\star)+h_{i_\star}(x_\star,l_\star,\beta_i)\Big\}\ge 0,~~\text{if}~~(x_\star,l_\star)\in(0,R)\times (0,K),\\
\partial_l\phi(0,l_\star)+\underset{\vartheta \in \mathcal{O}}{\inf} \Big\{\ds \sum_{i=1}^I\mathbb{S}_i(l_\star,\vartheta)\partial_x\phi_i(0,l_\star)+h_0(l_\star,\vartheta)\Big\} \leq 0,~~\text{if}~~x_\star=0,~~l_\star\in(0,K).
\end{cases}.
\end{align*}
b) We say that $v$ is a continuous viscosity sub solution of the $\mathcal{W}_{alsh}(\mathbb{S})$ system \eqref{eq system Walsh}, if for all test function $\phi \in \mathcal{C}^{2,0}_{{\bf 0},1}\big(\mathcal{N}_R\times [0,K]\big)$ and for all local maximum point 
$(x_\star,i_\star,l_\star)\in [0,R]\times [\![1,I]\!] \times [0,K]$ of $v-\phi$, with $(v-\phi)_{i_\star}(x_\star,l_\star)=0$, we have:
\begin{align*}
\begin{cases}
\lambda \phi_{i_\star}(x_\star,l_\star)+\underset{\beta_{i_\star}\in \mathcal{B}_{i_\star}}{\sup}\Big\{-\sigma_{i_\star}(x_\star,l_\star,\beta_i)\partial^2_x\phi_{i_\star}(x_\star,l_\star)+\\
b_{i_\star}(x_\star,l_\star,\beta_i)\partial_x\phi_{i_\star}(x_\star,l_\star)+h_{i_\star}(x_\star,l_\star,\beta_i)\Big\}\leq 0,~~\text{if}~~(x_\star,l_\star)\in(0,R)\times (0,K),\\
\partial_l\phi(0,l_\star)+\underset{\vartheta \in \mathcal{O}}{\inf} \Big\{\ds \sum_{i=1}^I\mathbb{S}_i(l_\star,\vartheta)\partial_x\phi_i(0,l_\star)+h_0(l_\star,\vartheta)\Big\} \ge 0,~~\text{if}~~x_\star=0,~~l_\star\in(0,K).
\end{cases}.
\end{align*}
c) We say that $u$ is a continuous viscosity solution of
the $\mathcal{W}_{alsh}(\mathbb{S})$ system \eqref{eq system Walsh}, if it is both a continuous viscosity super and sub solution of the $\mathcal{W}_{alsh}(\mathbb{S})$ system \eqref{eq system Walsh}.
\end{Definition}
The main result of this work is the following Theorem:
\begin{Theorem} (Comparison Theorem.)\label{th: compa theorem principale}
Assume assumption $(\mathcal{H})$. Let $v\in \mathcal{C}\big(\mathcal{N}_R\times [0,K]\big)$ a continuous viscosity sub solution and $u\in \mathcal{C}\big(\mathcal{N}_R\times [0,K]\big)$ a continuous viscosity super solution of the Walsh's spider  HJB system - $\mathcal{W}_{alsh}(\mathbb{S})$ - given in \eqref{eq system Walsh}, satisfying the following boundary conditions:
\begin{align*}
&\forall i\in [\![1,I]\!],~~\forall l\in[0,K],~~u_i(R,l)\ge v_i(R,l),
\\&\forall i\in [\![1,I]\!],~~\forall x\in[0,R],~~u_i(x,K)\ge v_i(x,K). 
\end{align*}
Then we have:
$$\forall (x,i,l)\in [0,R]\times [\![1,I]\!]\times [0,K],~~u_i(x,l)\ge v_i(x,l). $$
\end{Theorem}
We obtain also the version of Theorem of \ref{th: compa theorem principale}, in the case when the local time variable belongs to an unbounded domain. Its sketch of proof is given at the end of Section \ref{sec preuve theorem principal}.
\begin{Theorem}\label{th: compa theorem l infiny}
Assume now the dependency of the data $(\mathcal{D})$ are extended with respect to the variable $l$ to the unbounded domain $[0,+\infty)$, whereas assumption $(\mathcal{H})$ still holds true (replacing $K$ by $+\infty$). Let $v\in \mathcal{C}\big(\mathcal{N}_R\times [0,+\infty)\big)$ a continuous viscosity sub solution and $u\in \mathcal{C}\big(\mathcal{N}_R\times [0,+\infty)\big)$ a continuous viscosity super solution (in the sens of Definition \ref{def viscosity} with $K$ replaced by $+\infty$), of the following system:
\begin{eqnarray}\label{eq HJB bord infini}
&\begin{cases}
\textbf{HJB equation parameterized by the local-time on each ray:}\\ 
\lambda f_i(x,l)+\underset{\beta_i\in \mathcal{B}_i}{\sup}\Big\{-\sigma_i(x,l,\beta_i)\partial^2_xf_i(x,l)+
b_i(x,l,\beta_i)\partial_xf_i(x,l)+\\
h_i(x,l,\beta_i)\Big\}=0,~~(x,l)\in(0,R)\times (0,+\infty),\\
\textbf{Non linear local-time Kirchhoff's boundary condition at } \bf 0:\\
\partial_lf(0,l)+\underset{ \vartheta \in \mathcal{O}}{\inf} \Big\{\ds\sum_{i=1}^I\mathbb{S}_i(l,\vartheta)\partial_xf_i(0,l)+h_0(l,\vartheta)\Big\}=0,~~l\in(0,+\infty),\\
\textbf{Dirichlet boundary conditions outside } \bf 0:\\ 
f_i(R,l)=\chi_i(l),~~l\in[0,+\infty),\\
\textbf{Continuity condition at } \bf 0:\\ 
\forall (i,j)\in[\![1,I]\!]^2,~~\forall l\in [0,+\infty),~~f_i(0,l)=f_j(0,l).
\end{cases}
\end{eqnarray}
Assume that $u$ and $v$ satisfy a linear growth with respect to the variable $l$, namely there exists a constant $C>0$ such that:
$$\forall (x,l)\in[0,R]\times [0,+\infty),~~\forall i\in[\![1,I]\!],~~|u_i(x,l)|+|v_i(x,l)|\leq C(1+l).$$
Assume moreover that the following boundary condition is satisfied:
\begin{align*}
&\forall i\in [\![1,I]\!],~~\forall l\in[0,+\infty),~~u_i(R,l)\ge v_i(R,l). 
\end{align*}
Then we have:
$$\forall (x,i,l)\in [0,R]\times [\![1,I]\!]\times [0,+\infty),~~u_i(x,l)\ge v_i(x,l). $$   
\end{Theorem}
As a consequence of Theorem \ref{th: compa theorem l infiny}, we are able to obtain a comparison theorem in the case when the dependency with respect to the variable $l$ is removed, and this induces that a classical non linear Kirchhoff's boundary condition appears at $\bf 0$.
\begin{Theorem}\label{theo compa normal}
Assume assumption $(\mathcal{H})$ and that all the data $(\mathcal{D})$ have no dependence with respect to the variable $l$. 
Then the comparison theorem for continuous viscosity solutions (in the sens of Definition \ref{def viscosity} removing the dependency w.r.t to variable $l$) holds true for the following system posed on the star-shaped network $\mathcal{N}_R$:
\begin{eqnarray}\label{eq non linear HJB standard}
\begin{cases}
\textbf{HJB equation on each ray:}\\ 
\lambda u_i(x)+\underset{\beta_i\in \mathcal{B}_i}{\sup}\Big\{-\sigma_i(x,\beta_i)\partial^2_xu_i(x)+
b_i(x,\beta_i)\partial_xu_i(x)\\
+h_i(x,\beta_i)\Big\}=0,~~x\in(0,R),\\
\textbf{Non linear Kirchhoff's boundary condition at } \bf 0:\\
\underset{ \vartheta \in \mathcal{O}}{\inf} \Big\{\ds\sum_{i=1}^I\mathbb{S}_i(\vartheta)\partial_xu_i(0)+h_0(\vartheta)\Big\}=0,\\
\textbf{Dirichlet boundary conditions outside } \bf 0:\\ 
u_i(R)=\chi_i,\\
\textbf{Continuity condition at } \bf 0:\\ 
\forall (i,j)\in[\![1,I]\!]^2,~~u_i(0)=u_j(0).
\end{cases}
\end{eqnarray}
\end{Theorem}    
As we pointed out in the Introduction, recall that thanks to this new technique consisting of adding a variable representing the local time at $\bf 0$, we are able to obtain a comparison theorem for the system \eqref{eq non linear HJB standard} when it is not degenerate at $\bf 0$. Moreover, the Neumann (Kirchhoff) boundary condition is considered in the strong sense, without any dependency of the Hamiltonians at $\bf 0$, that is - up to our knowledge - a new advance for problems  with Neumann boundary conditions.\\
\\
\textbf{Proof of Theorem \ref{theo compa normal}}:
\begin{proof}
Let $v\in \mathcal{C}\big(\mathcal{N}_R\big)$ a continuous viscosity sub solution and $u\in \mathcal{C}\big(\mathcal{N}_R\big)$ a continuous viscosity super solution of \eqref{eq non linear HJB standard}, satisfying the following boundary conditions:
\begin{equation}\label{inega bord 1}
\forall i\in[\![1,I]\!],~~u_i(R)\ge v_i(R). 
\end{equation}
Assume by contradiction that:
\begin{eqnarray}\label{ref: contra classique normal}
 M=\sup \big\{~~v_i(x)-u_i(x),~~(x,i)\in \mathcal{N}_R~~\big\}>0.   
\end{eqnarray}
We will see that the last assumption will lead to a contradiction. As explained above, the central key is to use the comparison Theorem \ref{th: compa theorem l infiny}, having a local-time Kirchhoff's boundary condition at $\bf 0$.

We introduce for this purpose the following system now with non linear local-time Kirchhoff's boundary condition at $\bf 0$:
\begin{eqnarray}\label{eq non linear HJB inter}
\begin{cases}
\lambda f_i(x,l)+\underset{\beta_i\in \mathcal{B}_i}{\sup}\Big\{-\sigma_i(x,\beta_i)\partial^2_xf_i(x,l)+\\
b_i(x,\beta_i)\partial_xf_i(x,l)+
h_i(x,\beta_i)\Big\}=0,~~x\in(0,R)\times(0,+\infty),\\
\partial_lf_i(0,l)+\underset{ \vartheta \in \mathcal{O}}{\inf} \Big\{\ds\sum_{i=1}^I\mathbb{S}_i(\vartheta)\partial_xf_i(0,l)+h_0(\vartheta)\Big\}=0,~~l\in(0,+\infty)\\ 
f_i(R,l)=\chi_i,~~l\in[0,+\infty)\\
\forall (i,j)\in[\![1,I]\!]^2,~~f_i(0,l)=f_j(0,l),~~l\in[0,+\infty).
\end{cases}
\end{eqnarray}
We are going to show that the following maps:
\begin{align*}
&u^\varepsilon:=\begin{cases}
    \mathcal{N}_R\times [0,K]\to \R,\\
    \big((x,i),l\big)\mapsto u_i(x)+\varepsilon\exp(-\varepsilon l)
\end{cases}
&\\\text{and respectively,}\\
&v^\varepsilon:=\begin{cases}
    \mathcal{N}_R\times [0,K]\to \R,\\
    \big((x,i),l\big)\mapsto v_i(x)-\varepsilon\exp(-\varepsilon l)
\end{cases},
\end{align*}
are indeed super solution (resp. sub solution) of \eqref{eq non linear HJB inter}. Hence if the following boundary condition is satisfied:
\begin{align}\label{eq condition bords 1}
\forall i\in [\![1,I]\!],~~\forall l\in[0,+\infty),~~u_i^\varepsilon(R,l)\ge v_i^\varepsilon(R,l),
\end{align}
it follows from Theorem \ref{th: compa theorem l infiny} that:
\begin{align*}
&\forall i\in [\![1,I]\!],~~\forall (x,l)\in[0,R]\times[0,+\infty),~~u_i(x)+\varepsilon\exp(-\varepsilon l)\ge v_i(x)-\varepsilon\exp(-\varepsilon l). 
\end{align*}
Therefore sending $\varepsilon\searrow 0$ in the last equation, we obtain that:
\begin{align*}
&\forall i\in [\![1,I]\!],~~\forall x\in[0,R],~~u_i(x)\ge v_i(x), 
\end{align*}
and then a contradiction with \eqref{ref: contra classique normal}.
First it is easy to check that \eqref{eq condition bords 1} holds true, using \eqref{inega bord 1}. Moreover, both $u^\varepsilon$ and $v^\varepsilon$ are uniformly bounded in the domain $\mathcal{N}_R\times[0,+\infty)$ and satisfy then the growth assumption of Theorem \ref{th: compa theorem l infiny}. 
Let us now turn to prove that $u^\varepsilon$ is a super solution (resp. $v^\varepsilon$ is a sub solution) of \eqref{eq non linear HJB inter}. Let $\phi\in \mathcal{C}^{2,0}_{0,1}\big(\mathcal{N}_R\times[0,+\infty)\big)$ such that $u^\varepsilon-\phi$ has a local minimum point at $(x_\star,i_\star,l_\star)\in(0,R)\times [\![1,I]\!]\times (0,+\infty)$.\\
\textbf{Step 1} Assume that $x_\star>0$. There exists therefore an open set $\mathcal{V}$ of $(0,R)\times (0,+\infty)$ containing $(x_\star,l_\star)$ and strictly included in the ray $\mathcal{R}_{i_\star}$, such that:
$$\forall (x,l)\in\mathcal{V},~~u^\varepsilon_{i_\star}(x,l)-\phi_{i_\star}(x,l)\ge u^\varepsilon_{i_\star}(x_\star,l_\star)-\phi_{i_\star}(x_\star,l_\star).$$
In particular, for $l=l_\star$ we get that:
$$\forall x\in \mathcal{V}_{l_\star},~~u_{i_\star}(x)+\varepsilon\exp(-\varepsilon l_\star)-\phi_{i_\star}(x,l_\star)\ge u_{i_\star}(x_\star)+\varepsilon\exp(-\varepsilon l_\star)-\phi_{i_\star}(x_\star,l_\star),$$
where $\mathcal{V}_{l_\star}$ denotes the $l_\star$-level open set of $\mathcal{V}$, namely:
$$\mathcal{V}_{l_\star}:=\{x\in \mathcal{V},~~(x,l_\star)\in \mathcal{V}\}.$$
We conclude that $x\mapsto \phi_{i_\star}(x,l_\star)$ is a test function of $u$ at $(x_\star,i_\star)$ and since $u$ is a super solution of \eqref{eq non linear HJB standard}, we have:
$$\lambda \phi_{i_{\star}}(x_\star,l_\star)+\underset{\beta_{i_\star}\in \mathcal{B}_{i_{\star}}}{\sup}\Big\{-\sigma_{i_{\star}}(x_\star,\beta_{i_{\star}})\partial^2_x\phi_{i_{\star}}(x_\star,l_\star)+
b_{i_{\star}}(x_\star,\beta_{i_{\star}})\partial_x\phi_{i_{\star}}(x_\star,l_\star)+
h_{i_{\star}}(x_\star,\beta_{i_{\star}})\Big\}\ge 0.$$
\textbf{Step 2:} Assume that $x_\star=0$. There exists an open set $\mathcal{V}$ (for the geodesic metric of the network) containing $(0,l_\star)$, such that:
\begin{align}\label{eq Vamos 1}
\forall (x,l)\in\mathcal{V},~~\forall i\in[\![1,I]\!],~~u^\varepsilon_{i}(x,l)-\phi_{i}(x,l)\ge u^\varepsilon(0,l_\star)-\phi(0,l_\star).    
\end{align}
Using the special case when $l=l_\star$ we obtain with the same arguments used in \textbf{Step 1}, that:
\begin{align}\label{inegalite 1}
 \underset{ \vartheta \in \mathcal{O}}{\inf} \Big\{\ds\sum_{i=1}^I\mathbb{S}_i(\vartheta)\partial_x\phi_i(0,l_\star)+h_0(\vartheta)\Big\}\leq 0.   
\end{align}
Now if $x=0$, from \eqref{eq Vamos 1} we obtain:
\begin{align}\label{eq Vamos 2}
\forall l\in\mathcal{V}_0,~~u(0)+\varepsilon\exp(-\varepsilon l)-\phi(0,l)\ge u(0)+\varepsilon\exp(-\varepsilon l_\star)-\phi(0,l_\star),  
\end{align}
where the open set $\mathcal{V}_0$ is given by:
$$\mathcal{V}_0:=\{l\in \mathcal{V},~~(0,l)\in \mathcal{V}\}.$$
Hence, \eqref{eq Vamos 2} implies:
$$\partial_l\phi(0,l_\star)=\partial_l\big(\varepsilon\exp(-\varepsilon l)\big)_{l=l_\star}=-\varepsilon^2\exp(-\varepsilon l_\star)\leq 0.$$
It follows from \eqref{inegalite 1} that:
$$\partial_l\phi(0,l_\star)+\underset{ \vartheta \in \mathcal{O}}{\inf} \Big\{\ds\sum_{i=1}^I\mathbb{S}_i(\vartheta)\partial_x\phi_i(0,l_\star)+h_0(\vartheta)\Big\}\leq 0. $$
In conclusion $u^\varepsilon$ is super solution of \eqref{eq non linear HJB inter}. Same arguments lead to show that $v^\varepsilon$ is sub solution of \eqref{eq non linear HJB inter} and that achieves the proof. 
\end{proof}
\section{A short example for the construction of test functions in the linear case }\label{sec : exemple}
As a short introduction of the method that will be used in this work to prove our main Theorem \ref{th: compa theorem principale}, we propose in this Section for the convenience of the reader, a simple example in the linear framework. 

Consider the following elliptic linear PDE  posed on the open set $(0,R)$, with Neumann boundary condition at $x=0$ and Dirichlet boundary condition at $x=R$:
\begin{eqnarray}\label{eq EDP linear}
\begin{cases}
\lambda u(x)-\sigma(x)\partial^2_xu(x)=0,~~x\in(0,R),\\
\partial_xu(0)=0,~~u(R)=z,
\end{cases}
\end{eqnarray}
where $\lambda>0$, $\sigma\in \mathcal{C}[0,R]$ is strictly positive (elliptic) and $z\in \R$.
We are going to give a simple sketch of proof for a comparison theorem. 

Let $f$ (and resp. $g$) be a super (resp. sub) continuous viscosity solution of \eqref{eq EDP linear} (namely in the class $\mathcal{C}([0,R])$ in the sens of Definition \ref{def viscosity} adapted to the simple example \eqref{eq EDP linear}). As explained in Introduction, all our concentrations are focused on the behavior at the boundary $x=0$, hence we will only show that the following assumption:
\begin{eqnarray}\label{eq: contrad edp lin}
\sup_{x\in [0,R]} \big\{~g(x)-f(x)~\big\}=g(0)-f(0)>0,
\end{eqnarray}
will lead to a contradiction.

Let $\eta>0$ and $\gamma>0$ be two small parameters and $\varepsilon \in (0,R)$ designed to drive the construction of the test functions at the neighborhood $x=0$.
Define $\overline{\phi}=\overline{\phi}(\varepsilon,\eta,\gamma)$ the solution of the following ordinary differential equation:
\begin{eqnarray}
 \begin{cases}
 -\partial_x^2\overline{\phi}(x)+\ds \frac{\lambda f(x)}{\sigma(x)}=-\eta,~~x\in (0,\varepsilon),\\
\overline{\phi}(0)=0,~~\overline{\phi}(\varepsilon)=f(\varepsilon)-f(0)-\gamma.
\end{cases}
\end{eqnarray}
The solution satisfies:
$$\forall x\in [0,\varepsilon],~~\overline{\phi}(x)=f(\varepsilon)-f(0)-\gamma+\partial_x\overline{\phi}(0)(x-\varepsilon)+\int_\varepsilon^x\int_0^u\big(\eta+\frac{\lambda f(z)}{\sigma(z)}\big)dzdu.$$
We are going to prove that $\overline{\phi}$ is a test function of the super solution $f$ at $x=0$. To obtain this fact, it is enough to show that the minimum of $f-\overline{\phi}$ on the compact set $[0,\varepsilon]$ is necessary reached at $x=0$, and this is the case since:\\
-we have $f(0)-\overline{\phi}(0)<f(\varepsilon)-\overline{\phi}(\varepsilon),$\\
-and if the minimum of $f-\overline{\phi}$ is reached at the interior of $[0,\varepsilon]$, for instance at $y\in (0,\varepsilon)$, because $f$ is a super solution we should have:
$$ -\eta\sigma(y)=-\sigma(y)\partial_x^2\overline{\phi}(y)+\ds \lambda f(y)\ge 0,$$
which is a contradiction (recall that $\eta>0$ and $\sigma(y)>0$).
Therefore $\overline{\phi}$ is a test function of the super solution $f$ at $x=0$, and this implies:
$$\partial_x\overline{\phi}(0)\leq 0.$$
Using that $\overline{\phi}(0)=0$, we get:
$$0\ge \partial_x\overline{\phi}(0)\varepsilon =f(\varepsilon)-f(0)-\gamma+\int_\varepsilon^0\int_0^u\big(\eta+\frac{\lambda f(z)}{\sigma(z)}\big)dzdu,$$
and therefore:
\begin{eqnarray}\label{eq in lin 1}
\int_0^\varepsilon\int_0^u\big(\eta+\frac{\lambda f(z)}{\sigma(z)}\big)dzdu \ge f(\varepsilon)-f(0)-\gamma.
\end{eqnarray}
Analogously for the sub solution $g$, we define $\underline{\phi}=\underline{\phi}(\varepsilon,\eta,\gamma)$ the solution of the following ordinary differential equation:
\begin{eqnarray}
 \begin{cases}
 -\partial_x^2\underline{\phi}(x)+\ds \frac{\lambda g(x)}{\sigma(x)}=\eta,~~x\in (0,\varepsilon),\\
\underline{\phi}(0)=0,~~\underline{\phi}=g(\varepsilon)-g(0)+\gamma,
\end{cases}.   
\end{eqnarray}
The solution satisfies:
$$\forall x\in [0,\varepsilon],~~\underline{\phi}(x)=g(\varepsilon)-g(0)+\gamma+\underline{\phi}(x-\varepsilon)+\int_\varepsilon^x\int_0^u\big(-\eta+\frac{\lambda g(z)}{\sigma(z)}\big)dzdu,$$
and with the same arguments behind, we can show that $\underline{\phi}$ is a test function of the sub solution $g$ at $x=0$ (the maximum of $g-\underline{\phi}$ on the compact set $[0,\varepsilon]$ is necessary reached at $x=0$). We get then:
$$\partial_x\underline{\phi}(0)\ge 0,$$
and therefore using that $\underline{\phi}(0)=0$:
\begin{eqnarray}\label{eq in lin 2}
\int_0^\varepsilon\int_0^u\big(-\eta+\frac{\lambda g(z)}{\sigma(z)}\big)dzdu \leq g(\varepsilon)-g(0)+\gamma.
\end{eqnarray}
Now combining both \eqref{eq in lin 1} and \eqref{eq in lin 2}, we obtain:
$$
\int_0^\varepsilon\int_0^u\big(-2\eta+\frac{\lambda \big(g(z)-f(z))}{\sigma(z)}\big)dzdu \leq g(\varepsilon)-g(0)-f(\varepsilon)+f(0)+2\gamma.
$$
Sending $\gamma \searrow 0$ and $\eta \searrow  0$, we have:
$$
\limsup_{\eta \searrow  0}\int_0^\varepsilon\int_0^u\big(-2\eta+\frac{\lambda \big(g(z)-f(z))}{\sigma(z)}\big)dzdu \leq \liminf_{\gamma \searrow 0}\big(g(\varepsilon)-g(0)-f(\varepsilon)+f(0)+2\gamma\big),
$$
which leads to:
$$
\int_0^\varepsilon\int_0^u\frac{\lambda \big(g(z)-f(z))}{\sigma(z)}dzdu \leq g(\varepsilon)-g(0)-f(\varepsilon)+f(0).
$$
Therefore as expected at the beginning of the sketch of proof, if we assume \eqref{eq: contrad edp lin} that is:
\begin{eqnarray*}  
    \sup_{x\in [0,R]} \big\{~g(x)-f(x)~\big\}=g(0)-f(0)>0,
\end{eqnarray*}
we get:
$$
\int_0^\varepsilon\int_0^u\frac{\lambda \big(g(z)-f(z))}{\sigma(z)}dzdu \leq 0.
$$
Dividing by $\varepsilon^2$ we obtain:
$$
\frac{1}{\varepsilon^2}\int_0^\varepsilon\int_0^u\frac{\lambda \big(g(z)-f(z))}{\sigma(z)}dzdu \leq 0,
$$
and finally sending $\varepsilon \searrow  0$:
$$
\lim_{\varepsilon \searrow  0}\frac{1}{\varepsilon^2}\int_0^\varepsilon\int_0^u\frac{\lambda \big(g(z)-f(z))}{\sigma(z)}dzdu \leq 0,
$$
we get:
$$\frac{\lambda}{\sigma(0)}\big(g(0)-f(0)\big)\leq 0, $$
and this leads to a contradiction with the assumption \eqref{eq: contrad edp lin} using that $\lambda>0$ and the ellipticity condition at $x=0$: ($\sigma(0)>0$). We obtain therefore the result expected at the beginning of this sketch of proof.\\

\section{Proof of Theorem \ref{th: compa theorem principale}}\label{sec preuve theorem principal}
We state first the following Proposition, that will be useful for the construction of the test functions at the neighborhood of the vertex $\bf 0$.
\begin{Proposition}\label{pr: solva EDO paramatrique}
Let $r>0$ and $(B,H)\in \R^2$.
Fix $\varepsilon>0$, $\kappa>0$ and $(\eta,\gamma)\in \R^2$, be four small parameters $(\varepsilon,\kappa,|\eta|,|\gamma|<<1)$ (expected to be sent to 0.)
Assume that $\varepsilon$ is small enough in order to satisfy:
\begin{equation}\label{eq conditions petites 1}
1-|B|\varepsilon\exp(|B|\varepsilon)>0.
\end{equation}
Let $w^\kappa$ and $\big(z_i^{\varepsilon,\kappa}\big)_{i \in[\![1,I]\!]}$ be real bounded sequences indexed by $(\varepsilon,\kappa)$:
$$\exists w\ge0,~\exists z\ge 0,~~\sup_{\kappa \ge 0}|w^\kappa| \leq w,~~\underset{i \in[\![1,I]\!]}{\max}\sup_{\varepsilon\ge 0}\sup_{\kappa \ge 0}|z_i^{\varepsilon,\kappa}|\leq z.$$
For a given $\ell \in (0,+\infty)$ and $S\ge 0$:\\
(i)-the following parametric ordinary differential equation posed on the domain $\mathcal{N}_{\varepsilon}\times [\ell-\kappa,\ell+\kappa]$:
\begin{align}
\begin{cases}\label{eq EDO para}
r\psi_i(x,l)-\partial_x^2\psi_i(x,l)+B|\partial_x\psi_i(x,l)|+H+\eta=0,~~x\in (0,\varepsilon),~~l\in (\ell-\kappa,\ell+\kappa),\\
\psi(0,l)=w^{\kappa}+S(l-\ell),~~\psi_i(\varepsilon,l)=z_i^{\varepsilon,\kappa}+S(l-\ell)+\gamma,~~l\in [\ell-\kappa,\ell+\kappa],
\end{cases},
\end{align}
admits a unique solution:
$$\psi=\psi\Big(\varepsilon,\kappa,\eta,\gamma,S,w^\kappa,\big(z_i^{\varepsilon,\kappa}\big)_{i \in[\![1,I]\!]}\Big)$$
in the class $\mathcal{C}^{2,0}_{{\bf 0},1}\big(\mathcal{N}_{\varepsilon}\times [\ell-\kappa,\ell+\kappa]\big)$.\\
(ii)-Fix $\beta>0$. As soon as we impose $\kappa=\kappa_\varepsilon$ small enough in order to verify
\begin{align}\label{eq conditions petites 2}
1-\varepsilon\beta r\kappa\big(\exp(|B|\varepsilon)-1\big)-\varepsilon^2\beta\frac{r\kappa|B|\exp(2|B|\varepsilon)}{1-|B|\varepsilon\exp(|B|\varepsilon)}>0,
\end{align}
then there exists $$S(\beta)=S\Big(\beta,\varepsilon,\kappa,\eta,\gamma,w^\kappa,\big(z_i^{\varepsilon,\kappa}\big)_{i \in[\![1,I]\!]}\Big)\ge 0$$ 
such that:
\begin{equation}\label{eq condition aborption}
\forall l \in[\ell-\kappa,\ell+\kappa],~~\partial_l \psi(0,l)=S(\beta)\ge \varepsilon\beta\Big( |B||\partial_x\psi|+|H|+|\eta|\Big).
\end{equation}
(iii)-Assume moreover that the sequences $\Big(w^\kappa,\big(z_i^{\varepsilon,\kappa}\big)_{i \in[\![1,I]\!]}\Big)$ satisfy:
$$\forall i\in[\![1,I]\!],~~\lim_{\varepsilon \searrow 0}\limsup_{\kappa\searrow 0}\big|z_i^{\varepsilon,\kappa}-w^\kappa\big|=0.$$
Then for all $\beta>0$ and the choice of the parameter $S(\beta)=S\Big(\beta,\varepsilon,\kappa,\eta,\gamma,w^\kappa,\big(z_i^{\varepsilon,\kappa}\big)_{i \in[\![1,I]\!]}\Big)$ satisfying \eqref{eq condition aborption}, the solution:
$$\psi=\psi\Big(\varepsilon,\kappa,\eta,\gamma,S(\beta),w^\kappa,\big(z_i^{\varepsilon,\kappa}\big)_{i \in[\![1,I]\!]}\Big)=\psi^{\varepsilon,\kappa,\eta,\gamma}_\beta$$
of \eqref{eq EDO para} satisfies:
\begin{align}\label{eq clef finale}
\underset{\eta\searrow 0,\gamma \searrow 0}{\limsup_{\varepsilon \searrow 0,\kappa \searrow 0,}}\max_{i \in [\![1,I]\!]}\sup_{l\in [\ell-\kappa,\ell+\kappa]}\Big|\frac{2}{\varepsilon^2}\int_0^\varepsilon\int_{0}^u\psi_{i,\beta}^{\varepsilon,\kappa,\eta,\gamma}(z,l)dzdu-w^{\kappa}\Big|=0.
\end{align}
\begin{proof}
Fix for a while $l\in [\ell-\kappa,\ell+\kappa]$, viewed as an external parameter.
We obtain easily the unique solvability on each ray $\mathcal{R}_{\varepsilon,i}:=[0,\varepsilon]\times\{i\},i \in[\![1,I]\!]$ of \eqref{eq EDO para} (using for instance the Corollary 1.9-II  given in \cite{EDO}). Remarking that the sequence $w^\kappa$ is independent of $i \in[\![1,I]\!]$, we obtain by extension the unique solvability in the class $\mathcal{C}^{2}\big(\mathcal{N}_{\varepsilon}\big)$.
Set
\begin{equation}\label{bron phi uni}
M=M(\kappa,\gamma,\eta,S)=w+z+\kappa S+|\gamma|+\frac{|H|+|\eta|}{r}.   
\end{equation}
It is easy to check, using the comparison theorem for elliptic problems having Dirichlet boundary conditions, that the following constant map $x\mapsto M$ (resp. $x\mapsto- M$) is a super solution (resp. sub solution) of \eqref{eq EDO para}. We have then:
\begin{align}\label{eq born uni}
    \forall i \in[\![1,I]\!],~~ \forall (x,l)\in [0,\varepsilon]\times [\ell-\kappa,\ell+\kappa],~~|\psi_i(x,l)|\leq M.
\end{align}
Since we have for all $i \in[\![1,I]\!]$ and for all $(x,l)\in (0,\varepsilon)\times (\ell-\kappa,\ell+\kappa)$:
\begin{align*}
|\partial_x^2\psi_i(x,l)|\leq |B||\partial_x\psi_i(x,l)|+rM+|H|+|\eta|, 
\end{align*}
it follows from the classical differential version of Gr\"onwall's Lemma that for all $(x,l)\in (0,\varepsilon)\times (\ell-\kappa,\ell+\kappa)$:
\begin{align}\label{eq grad 0}
\nonumber &|\partial_x\psi_i(x,l)|\leq |\partial_x\psi_i(0,l)|\exp(\int_0^\varepsilon|B|dx)+\int_0^\varepsilon\big(r M+|H|+|\eta|\big)\exp(\int_x^\varepsilon|B|du)dx\\
&\leq |\partial_x\psi_i(0,l)|\exp(|B|\varepsilon)+\frac{r M+|H|+|\eta|}{|B|}\big(\exp(|B|\varepsilon)-1\big),    
\end{align}
(with the convention $\big(\big(\exp(|B|\varepsilon)-1\big)/|B|\big)_{|B|=0}=\varepsilon$).

Observe now that for all $l \in (\ell-\kappa,\ell+\kappa)$ and for all $i \in[\![1,I]\!]$:
\begin{eqnarray}\label{eq 1 deriv}
\gamma+z_i^{\varepsilon,\kappa}-w^\kappa =\varepsilon\partial_x\psi_i(0,l)+\int_0^\varepsilon\int_0^u\Big(r\psi_i(z,l)+B |\partial_x\psi_i(z,l)|+H+\eta\Big)dzdu,
\end{eqnarray}
and therefore from \eqref{bron phi uni}, \eqref{eq born uni} and \eqref{eq 1 deriv}:
\begin{eqnarray}\label{eq 2 deriv}
&|\partial_x\psi_i(0,l)|\leq \ds \frac{1}{\varepsilon}|\gamma+z_i^{\varepsilon,\kappa}-w^\kappa|+\varepsilon\exp(|B|\varepsilon)\Big(C_1(\kappa,\gamma,\eta,S)+|B||\partial_x\psi_i(0,l)|\Big),\\
\label{eq constant 1}
&C_1=C_1(\kappa,\gamma,\eta,S)=r(w+z+\kappa S+|\gamma|)+2(|H|+|\eta|).
\end{eqnarray}
From \eqref{eq 2 deriv},  we see that as soon as we impose $\varepsilon$ small enough to get:
$$1-|B|\varepsilon\exp(|B|\varepsilon)>0,$$
we obtain that for all $l\in (\ell-\kappa,\ell+\kappa)$ and for all $i\in[\![1,I]\!]$:
\begin{align}\label{eq grad en 0}
 |\partial_x\psi_i(0,l)|\leq \frac{1}{1-|B|\varepsilon\exp(|B|\varepsilon)}\Big[\frac{1}{\varepsilon}|\gamma+z_i^{\varepsilon,\kappa}-w^\kappa|+C_1\varepsilon\exp(|B|\varepsilon)\Big].   
\end{align}
Since the sequences $w^\kappa$ and $z_i^{\varepsilon,\kappa}$ are uniformly bounded and the parameter $\kappa<<1$ is small enough, we get from \eqref{eq EDO para}-\eqref{eq grad 0}-\eqref{eq grad en 0} that $\big( \partial^2_x\psi(\cdot,l)\big)_{l\in [\ell-\kappa,\ell+\kappa]}$ is uniformly bounded with respect to the parameter $l\in [\ell-\kappa,\ell+\kappa]$. Hence the sequences $\big( \psi(\cdot,l)\big)_{l\in [\ell-\kappa,\ell+\kappa]}$ and $\big( \partial_x\psi(\cdot,l)\big)_{l\in [\ell-\kappa,\ell+\kappa]}$ are Lipschitz equicontinuous uniformly with respect to the parameter $l\in [\ell-\kappa,\ell+\kappa]$. On the other hand, from \eqref{eq EDO para} and the triangle inequality, we obtain also that $\big( \partial^2_x\psi(\cdot,l)\big)_{l\in [\ell-\kappa,\ell+\kappa]}$ is Lipschitz equicontinuous uniformly with respect to the parameter $l\in [\ell-\kappa,\ell+\kappa]$. In other words:
\begin{equation}\label{eq Ascoli}\sup_{l\in [\ell-\kappa,\ell+\kappa]} \|\psi(\cdot,l)\|_{\mathcal{C}^{2+\alpha}(\mathcal{N}_\varepsilon)}\leq C_\varepsilon,\end{equation}
where $C_\varepsilon\ge 0$ is a positive constant independent of $\kappa$ and $\alpha \in (0,1)$.

Let us show that $\psi$ and its first-second derivative are continuous with respect to variable $l$. For this purpose, let $(l_n)$ be a sequence of $[\ell-\kappa,\ell+\kappa]$ converging to $l\in [\ell-\kappa,\ell+\kappa]$. We deduce with the aid of \eqref{eq Ascoli}, Ascoli's theorem and the boundary conditions, that $\psi(\cdot,l_n)$ will converge up to sub sequence in $\mathcal{C}^{2+\alpha}\big(\mathcal{N}_\varepsilon\big)$ to a solution of \eqref{eq EDO para}, that is indeed $\psi(\cdot,l)$ by uniqueness. Therefore the map:
$$\begin{cases}
 [\ell-\kappa,\ell+\kappa]\to \mathcal{C}^{2+\alpha}\big(\mathcal{N}_\varepsilon\big)\\
 l\mapsto \psi(\cdot,l),
\end{cases},$$
is continuous. We conclude that $\psi$ is in the class $\mathcal{C}^{2,0}_{{\bf 0},1}\big(\mathcal{N}_\varepsilon \times[\ell-\kappa,\ell+\kappa]\big)$, since clearly $l\mapsto \psi(0,l)\in \mathcal{C}^{1}\big([\ell-\kappa,\ell+\kappa]\big)$.

Let us show now \eqref{eq condition aborption}. Fix $\beta>0$. With the aid of \eqref{eq grad 0} and \eqref{eq grad en 0}, we have:
\begin{align}\label{eq: born grad finale  EDO} 
&\nonumber |\partial_x\psi|=\underset{i\in[\![1,I]\!]}{\max}\Big\{\sup\Big\{~~|\partial_x\psi_i(x,l)|,~~(x,l)\in[0,\varepsilon]\times[\ell-\kappa
,\ell+\kappa]~~\Big\}\Big\}\leq \\
&\nonumber  \frac{\exp(|B|\varepsilon)}{1-|B|\varepsilon\exp(|B|\varepsilon)}\Big[~~\underset{i\in[\![1,I]\!]}{\max}\frac{1}{\varepsilon}|\gamma+z_i^{\varepsilon,\kappa}-w^\kappa|+C_1(\kappa,\gamma,\eta,S)\varepsilon\exp(|B|\varepsilon)~~\Big] +\\
&\frac{C_1(\kappa,\gamma,\eta,S)}{|B|}\big(\exp(|B|\varepsilon)-1\big),
\end{align}
where we recall that $C_1=C_1(\kappa,\gamma,\eta,S)$ is given in \eqref{eq constant 1} by:
$$C_1(\kappa,\gamma,\eta,S)=r(w+z+\kappa S+|\gamma|)+2(|H|+|\eta|).$$
Let choose $\kappa_\varepsilon>0$ small enough, such that for all $\kappa\leq \kappa_\varepsilon$:
$$1-\varepsilon\beta r\kappa\big(\exp(|B|\varepsilon)-1\big)-\varepsilon^2\beta\frac{r\kappa|B|\exp(2|B|\varepsilon)}{1-|B|\varepsilon\exp(|B|\varepsilon)}>0.$$
We observe from \eqref{eq: born grad finale  EDO} and the expression of the constant $C_1(\kappa,\gamma,\eta,S)$, that if we set the parameter $S\ge 0$:
$$S=S(\beta)=S\Big(\beta,\varepsilon,\kappa,\eta,\gamma,w^\kappa,\big(z_i^{\varepsilon,\kappa}\big)_{i \in[\![1,I]\!]}\Big)\ge 0,$$
such that:
\begin{align}\label{eq expr absorption}
  &\nonumber S(\beta):=\ds \varepsilon\beta|B|\frac{ \ds \frac{\exp(|B|\varepsilon)}{1-|B|\varepsilon\exp(|B|\varepsilon)}\Big[~~\underset{i\in[\![1,I]\!]}{\max}\frac{1}{\varepsilon}|\gamma+z_i^{\varepsilon,\kappa}-w^\kappa|+C_2(\gamma,\eta)\varepsilon\exp(|B|\varepsilon)~~\Big] }{\ds1-\varepsilon\beta r\kappa\big(\exp(|B|\varepsilon)-1\big)-\varepsilon^2\beta\frac{r\kappa|B|\exp(2|B|\varepsilon)}{1-|B|\varepsilon\exp(|B|\varepsilon)}}\\
&+\varepsilon\beta\frac{\ds C_2(\gamma,\eta)\big(\exp(|B|\varepsilon)-1\big)+|H|+|\eta|}{\ds 1-\varepsilon\beta r\kappa\big(\exp(|B|\varepsilon)-1\big)-\varepsilon^2\beta\frac{r\kappa|B|\exp(2|B|\varepsilon)}{1-|B|\varepsilon\exp(|B|\varepsilon)}},
\end{align}
with:
\begin{eqnarray}\label{eq expr constaet 2}
C_2(\gamma,\eta)=r(w+z+|\gamma|)+2(|H|+|\eta|). 
\end{eqnarray}
we will obtain:
$$S(\beta)\ge \varepsilon \beta \big( |B||\partial_x\psi|+|H|+|\eta|\big),$$
namely \eqref{eq condition aborption} holds true.

To conclude observe first that since both sequences $w^\kappa$ and $z_i^{\varepsilon,\kappa}$ are uniformly bounded, we will obtain in \eqref{eq expr absorption} that the parameter:
$$S(\beta)=S\Big(\beta,\varepsilon,\kappa,\eta,\gamma,w^\kappa,\big(z_i^{\varepsilon,\kappa}\big)_{i \in[\![1,I]\!]}\Big)=S\big(\beta,\varepsilon,\kappa,\eta,\gamma\big)$$
will satisfy:
\begin{equation}\label{eq convergence clef}
 \underset{\gamma \searrow 0}{\limsup_{\kappa \searrow 0,\eta \searrow 0}}~~\kappa S\big(\beta,\varepsilon,\kappa,\eta,\gamma\big)=\limsup_{\kappa \searrow 0}\limsup_{\eta \searrow 0}\limsup_{\gamma \searrow 0}~~\kappa S\big(\beta,\varepsilon,\kappa,\eta,\gamma\big)=0.  
\end{equation}
Let $l\in [\ell-\kappa,\ell+\kappa]$. Now for $\psi=\psi_\beta^{\varepsilon,\kappa,\eta,\gamma}$, which satisfies:
$$\partial_l\psi(0,l)=S\big(\beta,\varepsilon,\kappa,\eta,\gamma\big),$$
we have for all $i\in [\![1,I]\!]$ (recall that: $\psi_{\beta}^{\varepsilon,\kappa,\eta,\gamma}(0,l)=w^\kappa+S\big(\beta,\varepsilon,\kappa,\eta,\gamma\big)(l-\ell)$):
\begin{align*}
&\big|\frac{2}{\varepsilon^2}\int_0^\varepsilon\int_{0}^u\psi_{i,\beta}^{\varepsilon,\kappa,\eta,\gamma}(z,l)dzdu-w^\kappa\big|\leq\\
&\big|\frac{2}{\varepsilon^2}\int_0^\varepsilon\int_0^u\Big(\psi_{i,\beta}^{\varepsilon,\kappa,\eta,\gamma}(z,l)-\big(w^\kappa+S\big(\beta,\varepsilon,\kappa,\eta,\gamma\big)(l-\ell)\big)\Big)dzdu\big|+\kappa S\big(\beta,\varepsilon,\kappa,\eta,\gamma\big)\leq \\
&\big|\frac{2}{\varepsilon^2}\int_0^\varepsilon\int_0^u\int_0^z\partial_x\psi_{i,\beta}^{\varepsilon,\kappa,\eta,\gamma}(t)dtdzdu\big|+\kappa S\big(\beta,\varepsilon,\kappa,\eta,\gamma\big)\leq \\
&2\varepsilon|\partial_x\psi_\beta^{\varepsilon,\kappa,\eta,\gamma}|+\kappa S\big(\beta,\varepsilon,\kappa,\eta,\gamma\big).
\end{align*}
Hence if the sequences $\Big(w^\kappa,\big(z_i^{\varepsilon,\kappa}\big)_{i \in[\![1,I]\!]}\Big)$ satisfy:
$$\forall i\in[\![1,I]\!],~~\lim_{\varepsilon \searrow 0}\limsup_{\kappa\searrow 0}\big|z_i^{\varepsilon,\kappa}-w^\kappa\big|=0,$$
we see from the expression given \eqref{eq: born grad finale  EDO} and the convergence \eqref{eq convergence clef}, that we will obtain:
\begin{align*}
\underset{\eta\searrow 0,\gamma \searrow 0}{\limsup_{\varepsilon \searrow 0,\kappa \searrow 0,}}\max_{i \in [\![1,I]\!]}\sup_{l\in [\ell-\kappa,\ell+\kappa]}\Big|\frac{2}{\varepsilon^2}\int_0^\varepsilon\int_{0}^u\psi_{i,\beta}^{\varepsilon,\kappa,\eta,\gamma}(z,l)dzdu-w^{\kappa}\Big|=0.
\end{align*}
The proof is complete.
\end{proof}
\end{Proposition}
We are able now to prove the central result of this work: Theorem \ref{th: compa theorem principale}.
\begin{proof}
Let $f$ (and resp. $g$) be a super (resp. sub) continuous viscosity solution of \eqref{eq EDP linear}. \textbf{Fix in the sequel $\underline{\ell}\in [0,K)$.} We argue by contradiction assuming that:
 $$\sup\Big\{~g_i(x)-f_i(x),~~\big((x,i),l\big)\in \mathcal{N}_R\times[\underline{\ell},K]~\Big\}>0.$$
 Since $f$ and $g$ are in the class $\mathcal{C}\big(\mathcal{N}_R\times[0,K]\big)$, using the boundary condition given in the assumptions of the Theorem, the last supremum is necessary reached at a point:
 $$\big(x_\star,i_\star,l_\star\big)\in [0,R)\times[\![1,I]\!]\times [\underline{\ell},K).$$
\textbf{Step 1:  Classical arguments at the interior of each ray.} Assume first that: $x_\star>0$. Let $\mathcal{V}(x_\star)$ be a neighborhood of $x_\star$ strictly included in the ray $\mathcal{R}_{i\star}$. When $l_\star\in [\underline{\ell},K)$ is fixed, the following Hamiltonian (parameterized by $l_\star$) :
$$H_{i_\star}^{l_\star}:=\begin{cases}
\overline{\mathcal{V}(x_\star)}\times \R^3\to \R\\
 (x,u,p.S)\mapsto \lambda u+\sup_{\beta_{i_\star}\in B_{i_\star}}\Big\{-\sigma_{i_\star}^2(x,l_\star,\beta_{i_\star})S+\\
 b_{i_\star}(x,l_\star,\beta_{i_\star})p+h_{i_\star}(x,l_\star,\beta_{i_\star})\Big\},   
\end{cases},$$
is continuous because:\\
-linear in $u$,\\
-convex in $(p,S)$,\\
-assumption $(\mathcal{H}-\bf R)$ states that all the coefficients $(\sigma_i,b_i,h_i)_{i\in[\![1,I]\!]}$ are Lipschitz continuous, uniformly in the control variables $(\beta_i\in \mathcal{B}_i)_{i\in[\![1,I]\!]}$. This imply the Lipschitz continuity of $H_{i_\star}^{l_\star}$ with respect to the variable $x$.

Moreover, we have that the classical assumptions introduced in the seminal work \cite{User guide} (Theorem 3.3)  hold true, which are:
\begin{align*}
 &(i)~~\forall (x,u,v,p,S)\in [0,R]\times \R^3,~~\text{if}~u\ge v,~~\text{then}
 \\&H_{i_\star}^{l_\star}(x,u,p,S)-H_{i_\star}^{l_\star}(x,v,p.S)\ge \lambda(u-v),~~\text{with}~~\lambda >0,  
\end{align*}
\begin{align*}
 &(ii)~~\exists \omega \in  \mathcal{C}(\R_+,\R_+),~~\omega(0)=0,~~
 \forall \alpha>0,~~\forall (x,y,u,p,X,Y)\in [0,R]^2\times \R^4,\\
 &~\text{such that:}~~-3\alpha \begin{pmatrix}1&0\\
0&1
 \end{pmatrix}\leq\begin{pmatrix}X&0\\
 0&-Y
 \end{pmatrix}\leq 3\alpha\begin{pmatrix}1&-1\\
-1&1
 \end{pmatrix},~~\text{then}:\\
&H_{i_\star}^{l_\star}(y,u,\alpha(x-y).Y)-H_{i_\star}^{l_\star}(x,u,\alpha(x-y),X)\leq \omega\big(\alpha|x-y|^2+|x-y|\big).
\end{align*}
(Recall that the last property $(ii)$ holds true since the coefficients $(\sigma_i,b_i,h_i)_{i\in[\![1,I]\!]}$, are Lipschitz continuous uniformly in the control variables $(\beta_i\in \mathcal{B}_i)_{i\in[\![1,I]\!]}$ - see Example 3.6 in \cite{User guide} for instance - and $l=l_\star$ is fixed).

Since the Hamiltonian $H_{i_\star}^{l_\star}$ does not have any dependency with some derivative with respect to the variable $l$, we can proceed with the classical arguments of the theory of viscosity solution to obtain a contradiction.
Let us detail a little more this point since our class of test function is larger then the class $\mathcal{C}^2\big([0,R]\big)$ used in the classical framework. Indeed if $\phi\in \mathcal{C}^{0,2}_{0,1}\big(\mathcal{N}_R\times[0,K]\big)$ is a test function such that $f-\phi$ has a local minimum point at $(x_\star,i_\star,l_\star)\in(0,R)\times [\![1,I]\!]\times (0,K)$, there exists therefore an open set $\mathcal{V}$ of $(0,R)\times (0,K)$ containing $(x_\star,l_\star)$ and strictly included in the ray $\mathcal{R}_{i_\star}$, such that:
$$\forall (x,l)\in\mathcal{V},~~f_{i_\star}(x,l)-\phi_{i_\star}(x,l)\ge f_{i_\star}(x_\star,l_\star)-\phi_{i_\star}(x_\star,l_\star).$$
In particular, for $l=l_\star$ we get that:
$$\forall x\in \mathcal{V}_{l_\star},~~f_{i_\star}(x,l_\star)-\phi_{i_\star}(x,l_\star)\ge f_{i_\star}(x_\star,l_\star)-\phi_{i_\star}(x_\star,l_\star),$$
where $\mathcal{V}_{l_\star}$ denotes the $l_\star$-level open set of $\mathcal{V}$, defined by:
$$\mathcal{V}_{l_\star}:=\{x\in \mathcal{V},~~(x,l_\star)\in \mathcal{V}\}.$$
On the other hand by definition:
\begin{eqnarray*}
 &\lambda \phi_{i_\star}(x_\star,l_\star)+\underset{\beta_{i_\star}\in \mathcal{B}_{i_\star}}{\sup}\Big\{-\sigma_{i_\star}(x_\star,l_\star,\beta_{i_\star})\partial^2_x\phi_{i_\star}(x_\star,l_\star)+\\
&b_{i_\star}(x_\star,l_\star,\beta_{i_\star})\partial_x\phi_{i_\star}(x_\star,l_\star)+h_{i_\star}(x_\star,l_\star,\beta_{i_\star})\Big\}\ge 0.  
\end{eqnarray*}
This implies that that $x\mapsto \phi_{i_\star}(x,l_\star)$ is a test function of $x\mapsto f_{i_\star}(x,l_\star)$ in the neighborhood of $x_\star$, and $x\mapsto f_{i_\star}(x,l_\star)$ is a also a super solution of: 
$$H_{i_\star}^{l_\star}(x,u(x),\partial_xu(x),\partial_x^2u(x))=0,$$
at the interior of the ray $\mathcal{R}_{i_\star}$. The same arguments hold true for the sub solution $g$.
Hence we can proceed like in the classical case:\\
-applying the doubling variable method with the function (here parameterized by $(i_\star.l_\star)$) defined by :
$$\forall \varepsilon>0,~~w_\varepsilon^{(i_\star.l_\star)}(x,y)=g_{i_\star}(x,l_\star)-f_{i_\star}(x,l_\star)-\ds \frac{1}{2\varepsilon^2}|x-y|^2,~~(x,y)\in \mathcal{V}_{i_\star}(x_\star)\times\mathcal{V}_{i_\star}(x_\star),$$
(recall that $\mathcal{V}_{i_\star}(x_\star)$ is a neighborhood of $x_\star$ strictly included in $\mathcal{R}_{i_\star}$).\\
-use the Ishii's matrix lemma (see for example Theorem 3.2 in \cite{User guide}). (Recall that $H_{i_{\star}}^{l_\star}$ is continuous at the neighborhood of $x_\star$, which implies that the equivalent definition of a super (resp. sub solutions) with the closure of the second-order {\it super jet} (resp. {\it sub jet}) of $f$ (resp. $g$) at $x_\star$ holds true).\\
\textbf{In the rest of this proof we assume that:}
\begin{eqnarray}\label{eq: contrad visco en 0}
   &\sup \Big\{~g_i(x,l)-f_i(x,l),~\big((x,i),l\big)\in \mathcal{N}_R\times [\underline{\ell},K]~\Big\}=g(0,l_\star)-f(0,l_\star)>0,\\
   \nonumber &\text{where we recall that:}~~l_\star\in [\underline{\ell},K),~~\underline{\ell}\in(0,K).
\end{eqnarray}
\textbf{Step 2: Introduction of test functions depending on the speed of the Hamiltonians.} We scale first $f$ and $g$ at the vertex $\bf0$. Set: 
$$\Theta(f,g)=\frac{1}{2}\big(f(0,l_\star)+g(0,l_\star)\big).$$
It is easy to verify that:
\begin{eqnarray}\label{eq fonction scaller}
 u=f-\Theta(f,g),~~v=g-\Theta(f,g),  
\end{eqnarray}
are respectively super solution and sub solution of the following system with {\it non linear local time's Kirchhoff's boundary transmission}, posed on the domain $\mathcal{N}_R\times[0,K]$:
\begin{eqnarray}\label{eq: HJB scaller}
\begin{cases}
\lambda w_i(x,l)+\lambda\Theta(f,g)+\underset{\beta_i\in \mathcal{B}_i}{\sup}\Big\{-\sigma_i(x,l,\beta_i)\partial^2_xw_i(x,l)+\\
b_i(x,l,\beta_i)\partial_xw_i(x,l)+h_i(x,l,\beta_i)\Big\}=0,~~(x,l)\in(0,R)\times (0,K),\\
\partial_lw(0,l)+\underset{\vartheta \in \mathcal{O}}{\inf} \Big\{\ds\sum_{i=1}^I\mathbb{S}_i(l,\vartheta)\partial_xw_i(l,0)+h_0(l,\vartheta)\Big\}=0,~~l\in(0,K)\\
w_i(R,l)=\chi_i(l)-\Theta(f,g),~~l\in[0,K],\\
w_i(x,K)=\mathcal{T}_i(x)-\Theta(f,g),~~x\in[0,R]\\
\forall (i,j)\in[\![1,I]\!]^2,~~\forall l\in [0,K],~~w_i(0,l)=w_j(0,l).
\end{cases}~~
\end{eqnarray}
Remark that:
\begin{align}\label{eq positiv u v}
\nonumber &u(0,l_\star)=\frac{1}{2}\big(f(0,l_\star)-g(0,l_\star)\big)<0,~~v(0,l_\star)=\frac{1}{2}\big(g(0,l_\star)-f(0,l_\star)\big)>0,\\
&v(0,l_\star)-u(0,l_\star)=g(0,l_\star)-f(0,l_\star)>0.
\end{align}

In the sequel, drawing on the method introduced in the last Section \ref{sec : exemple}, we will build in the neighborhood of $(0,l_\star)$ test functions for the super solution $u$ (resp. sub solution $v$) of \eqref{eq: HJB scaller} that will solve ordinary differential equations possessing only constant coefficients. This coefficients may be viewed
as a kind a supreme envelope of all possible first order errors of the speed of the Hamiltonians $\Big(\mathcal{S}_{peed}(H_i)_{i\in[\![1,I]\!]}\!]\Big)$ defined by:
\begin{align*}
\mathcal{S}_{peed}(H_i):=\begin{cases}
    [0,R]\times[0,K]\times \R\to \R,\\
    (x,l,p)\mapsto \underset{\beta_i\in \mathcal{B}_i}{\sup}\Big\{\ds \frac{b_i(x,l,\beta_i)p+h_i(x,l,\beta_i)}{\sigma_i(x,l,\beta_i)}\Big\}
\end{cases},
\end{align*}
where the Hamiltonian $\Big(H_i,~i\in[\![1,I]\!]\Big)$ are given by:
\begin{align*}
H_i:=\begin{cases}
    [0,R]\times[0,K]\times \R^3\to \R,\\
    (x,l,u,p,S)\mapsto \lambda u+\underset{\beta_i\in \mathcal{B}_i}{\sup}\Big\{\ds-\sigma_i(x,l,\beta_i)S+b_i(x,l,\beta_i)p+h_i(x,l,\beta_i)\Big\}.
\end{cases}   
\end{align*}
The key point in the construction is to impose a derivative with respect to the 'local time' variable $l$ at $x=0$ that will absorb all the errors induced by {\it the Kirchhoff's speed of the Hamiltonians}, given by:
\begin{align}\label{eq Kir speed}
\mathcal{K}_{f}\mathcal{S}_{peed}\big((H_i)_{i\in[\![1,I]\!]}\big):=\begin{cases}
    [0,R]\times[0,K]\times \R\to \R,\\
    (x,l,p)\mapsto \ds \inf_{\vartheta \in \mathcal{O}}\Big\{\sum_{i=1}^I\mathbb{S}_i(l,\vartheta)\underset{\beta_i\in \mathcal{B}_i}{\sup}\Big\{\ds \frac{b_i(x,l,\beta_i)p+h_i(x,l,\beta_i)}{\sigma_i(x,l,\beta_i)}\Big\}\Big\}
\end{cases}.
\end{align}
Fix $\varepsilon>0$ and $\kappa>0$ two small parameters expected to be sent to $0$. Using \eqref{eq positiv u v} and the continuity of $u$ and $v$, there exists a neighborhood of the vertex $(0,l_\star)$ denoted by $\mathcal{V}\Big((0,l_\star),(\varepsilon,\kappa)\Big)$ and defined by:
$$\mathcal{V}\Big((0,l_\star),(\varepsilon.\kappa) \Big):=\Big\{\big((x,i),l\big)\in \mathcal{N}_R\times [\underline{\ell},K],~~x\leq \varepsilon,~~0<l_\star-\kappa\leq l \leq  l_\star+\kappa<K\Big\},$$
such that:
\begin{align}\label{eq positi u et v}
\forall\big((x,i),l\big)\in \mathcal{V}\Big((0,l_\star),(\varepsilon,\kappa)\Big),~~v_i(x,l)\ge 0,~~u_i(x,l)\leq 0,~~v_i(x,l)-u_i(x,l)\ge 0.
\end{align}
In the sequel we introduce also $\eta>0$ and $\gamma>0$ two small strictly positive parameters, designed to drive the construction of the test functions at the neighborhood of the vertex $(0,l_\star)$.

Set in the sequel:
\begin{align}\label{eq  expression const scal}
\nonumber \overline{t}\big(\lambda,\Theta(f,g)\big)=\lambda\Theta(f,g)\big(\frac{1}{\underline{\sigma}}\mathbf{1}_{\Theta(f,g)>0}+\frac{1}{\overline{\sigma}}\mathbf{1}_{\Theta(f,g)\leq 0}\big),\\
\underline{t}\big(\lambda,\Theta(f,g)\big)=\lambda\Theta(f,g)\big(\frac{1}{\underline{\sigma}}\mathbf{1}_{\Theta(f,g)\leq 0}+\frac{1}{\overline{\sigma}}\mathbf{1}_{\Theta(f,g)> 0}\big),\\
\nonumber \overline{u}^\kappa(0)=\sup\big\{u(0,l),~~l\in[l_\star-\kappa,l_\star+\kappa]\big\},\\
\nonumber\forall i\in [\![1,I]\!],~~
\underline{u}_i^\kappa(\varepsilon)=\inf\big\{u_i(\varepsilon,l),~~l\in[l_\star-\kappa,l_\star+\kappa]\big\},\\
\nonumber \underline{v}^\kappa(0)=\inf\big\{v(0,l),~~l\in[l_\star-\kappa,l_\star+\kappa]\big\},\\  \forall i\in [\![1,I]\!],~~
\overline{v}_i^\kappa(\varepsilon)=\sup\big\{v_i(\varepsilon,l),~~l\in[l_\star-\kappa,l_\star+\kappa]\big\},
\end{align}
Let $\overline{S}\ge 0$ and $\underline{S}\ge 0$ be two parameters designed to characterize the 'local time' derivative of the tests functions at $\bf 0$.
Using Proposition \ref{pr: solva EDO paramatrique}, we introduce \begin{align*}
&\overline{\phi}=\overline{\phi}(u,\varepsilon,\kappa,\eta,\gamma,\overline{S}),~~\underline{\phi}=\underline{\phi}(v,\varepsilon,\kappa,\eta,\gamma,\underline{S})
\end{align*}
(denoted in the next lines lines $\big(\overline{\phi},\underline{\phi}\big)$ for the seek of clarity) the two solutions of the  two following ordinary parametric differential equation systems posed on $\mathcal{N}_\varepsilon\times [l_\star-\kappa,l_\star+\kappa]$
\begin{align}\label{expression EDO sur sol}
\begin{cases}
\ds \frac{\lambda }{\overline{\sigma}}\overline{\phi}_i-\partial_x^2\overline{\phi}_i(x,l)+\overline{t}\big(\lambda,\Theta(f,g)\big)\\
 \ds \hspace{2,0 cm}+\Big(|b||\partial_x\overline{\phi}_i(x,l)|+|h|\Big)\big/\underline{\sigma}=-\eta,~~(x,l)\in (0,\varepsilon)\times(l_\star-\kappa,l_\star+\kappa),\\
\overline{\phi}(0,l)=\overline{u}^\kappa(0)+\overline{S}(l-l_\star),\\
\overline{\phi}_i(\varepsilon,l)=\underline{u}_i^\kappa(\varepsilon)-\gamma+\overline{S}(l-l_\star),~~l\in[l_\star-\kappa,l_\star+\kappa],~~i\in [\![1,I]\!],
\end{cases}
\end{align}
and:
\begin{align} \label{expression EDO sous sol}
\begin{cases}
\ds \frac{\lambda }{\overline{\sigma}}\underline{\phi}_i-\partial_x^2\underline{\phi}_i(x,l)+\underline{t}\big(\lambda,\Theta(f,g)\big)\\
 \ds \hspace{2,0 cm}-\Big(|b||\partial_x\underline{\phi}_i(x,l)|+|h|\Big)\big/\underline{\sigma}=\eta,~~(x,l)\in (0,\varepsilon)\times(l_\star-\kappa,l_\star+\kappa),\\
\underline{\phi}(0,l)=\underline{v}^\kappa(0)-\underline{S}(l-l_\star),\\
\underline{\phi}_i(\varepsilon,l)=\overline{v}_i^\kappa(\varepsilon)+\gamma-\underline{S}(l-l_\star),~~l\in[l_\star-\kappa,l_\star+\kappa],~~i\in [\![1,I]\!].
\end{cases}
\end{align}
Recall that $\underline{\sigma}$ and $\overline{\sigma}$ are given in Assumption $(\mathcal{H})$: $({\bf E})$ -$({\bf R}-ii)$.
Assumption $(\mathcal{H})$ combined with Proposition \ref{pr: solva EDO paramatrique} state that both solutions $\overline{\phi}$ and $\underline{\phi}$ are unique and in the class of test functions $\mathcal{C}^{2,0}_{{\bf  0},1}\big(\mathcal{N}_\varepsilon\times [l_\star-\kappa,l_\star+\kappa]\big)$.

\textbf{Step 3: $\overline{\phi}$ (resp. $\underline{\phi}$) is a test function of the super solution $u$ (resp. the sub solution $v$) at the vertex $\bf 0$ of the Walsh's spider HJB system \eqref{eq: HJB scaller}.}

We start to show that $\overline{\phi}$ is a test function of the super solution $u$  at $\bf 0$ of the HJB system \eqref{eq: HJB scaller}. The second case involving $\underline{\phi}$ for the sub solution $v$ can be treated with the same arguments.

We are then going to show that the minimum of $u-\overline{\phi}$ on the domain $\mathcal{N}_\varepsilon\times [l_\star-\kappa, l_\star+\kappa]$, is necessarily reached at some point $(0,\overline{l}_{\kappa})$ with $\overline{l}_{\kappa}\in [l_\star-\kappa,l_\star+\kappa]$. Using the boundary conditions satisfied by $\overline{\phi}$ and the expressions given in \eqref{eq  expression const scal}, we obtain that $\forall i \in[\![1,I]\!]$ and $\forall l \in [l_\star-\kappa, l_\star+\kappa]$:
\begin{align*}
 &u(0,l)-\overline{\phi}(0,l)=u(0,l)-\overline{u}^{\kappa}(0)-\overline{S}(l-l_\star)\leq -\overline{S}(l-l_\star),\\
&\forall i \in[\![1,I]\!],~~u_i(\varepsilon,l)-\overline{\phi}_i(\varepsilon,l)=u_i(\varepsilon,l)-\underline{u}^{\kappa}_i(\varepsilon)-\overline{S}(l-l_\star)+\gamma>-\overline{S}(l-l_\star),  \end{align*}
which implies that the minimum of $u-\overline{\phi}$ can not be reached at $(x=\varepsilon,l)$ with $l\in [l_\star-\kappa,l_\star+\kappa]$.

Assume now that minimum of $u-\overline{\phi}$ is reached at some point $$\big(y,j,\ell\big)\in (0,\varepsilon)\times [\![1,I]\!]\times [l_\star-\kappa,l_\star+\kappa].$$
As it is classical in the viscosity formulation, we can without lose of generality assume that:
$$u_j(y,\ell)=\overline{\phi}_{j}(y,\ell).$$
Since $u$ is a super solution  of the system \eqref{eq: HJB scaller}, we have by definition (recall that $0<l_\star-\kappa<l_\star+\kappa<K$):
\begin{align}\label{eq premier jet}
&\nonumber\lambda u_j(y,\ell)+\lambda \Theta(f,g)+\\
&\sup_{\beta_j\in \mathcal{B}_j}\Big\{-\sigma_j(y,\ell,\beta_j)\partial_x^2\overline{\phi}_{j}(y,\ell)+b_j(y,\ell,\beta_j)\partial_x\overline{\phi}_{j}(y,\ell)+h_j(y,\ell,\beta_j)\Big\}\ge 0. 
\end{align}
Hence using the expression of the system of ordinary differential equation \eqref{expression EDO sur sol} satisfied by $\overline{\phi}$  we get:
\begin{align}\label{eq expr contradiction}
&\nonumber \lambda u_j(y,\ell)+\lambda \Theta(f,g)+\sup_{\beta_j\in \mathcal{B}_j}\Big\{-\sigma_j(y,\ell,\beta_j)\eta- \sigma_j(y,\ell,\beta_j)\frac{\lambda u_j(y,\ell)}{\overline{\sigma}}\\&\nonumber -\sigma_j(y,\ell,\beta_j)\overline{t}\big(\lambda, \Theta(f,g)\big)-\sigma_j(y,\ell,\beta_j)\big(|b||\partial_x\overline{\phi}_{j}(y,\ell)|+|h|\big)\big/\underline{\sigma}\\&+b_j(y,\ell,\beta_j)\partial_x\overline{\phi}_{j}(y,\ell)+h_j(y,\ell,\beta_j)\Big\}\ge 0.  
\end{align}
Using now that $\eta>0$, $u_j(y,\ell)\leq 0$ (see \eqref{eq positi u et v}), the expression of $\overline{t}\big(\lambda,\Theta(f,g)\big)$ given in \eqref{eq  expression const scal} and the central ellipticity condition satisfied by the coefficient $\sigma_j$ (assumption $\mathcal{H}-\bf E )$, we get that the last quantity \eqref{eq expr contradiction} is smaller than:
\begin{align*}
&\Big[\lambda u_j(y,\ell)+\sup_{\beta_j\in \mathcal{B}_j}\{- \sigma_j(y,\ell,\beta_j)\frac{\lambda u_j(y,\ell)}{\overline{\sigma}}\}\Big]\\
&+\Big[\lambda \Theta(f,g)+\sup_{\beta_j\in \mathcal{B}_j}\{-\sigma_j(y,\ell,\beta_j)\overline{t}\big(\lambda, \Theta(f,g)\big)\}\Big]\\
&+\Big[\big(|b||\partial_x\overline{\phi}_{j}(y,\ell)|+|h|\big)-\big(|b||\partial_x\overline{\phi}_{j}(y,\ell)|+|h|\big)\Big]\\
&+\sup_{\beta_j\in \mathcal{B}_j}\{-\sigma_j(y,\ell,\beta_j)\eta\}\leq- \underline{\sigma}\eta<0, 
\end{align*}
and that leads to a contradiction with \eqref{eq premier jet}. 
We conclude that $\overline{\phi}$ is a test function of the super solution $u$ at the vertex $\bf 0$, namely there exists $\overline{l}_\kappa\in [l_\star-\kappa,l_\star+\kappa]$ such that:
\begin{eqnarray}\label{eq second jet u}
\partial_l\overline{\phi}(0,\overline{l}_\kappa)+\inf_{\vartheta \in \mathcal{O}}\Big\{\ds \sum_{i=1}^I\mathbb{S}_i(\overline{l}_\kappa,\vartheta)\partial_x\overline{\phi}_i(0,\overline{l}_\kappa)+h_0(\overline{l}_\kappa,\vartheta)\Big\}\leq 0.  
\end{eqnarray}
As announced at the beginning of this Step, using the same arguments, we can show easily that $\underline{\phi}$ is a test function of the sub solution $v$ at $\bf 0$, namely there exists $\underline{l}_\kappa\in [l_\star-\kappa,l_\star+\kappa]$ such that:
\begin{eqnarray}\label{eq second jet v}
\partial_l\underline{\phi}(0,\underline{l}_\kappa)+\inf_{\vartheta \in \mathcal{O}}\Big\{\ds \sum_{i=1}^I\mathbb{S}_i(\underline{l}_\kappa,\vartheta)\partial_x\underline{\phi}_i(0,\underline{l}_\kappa)+h_0(\underline{l}_\kappa,\vartheta)\Big\}\ge 0.  
\end{eqnarray}
\textbf{Step 3: Conclusion.} Fix $\vartheta\in \mathcal{O}$. Without loss of generality, in the construction of the tests functions in \textbf{Step 1} and the results obtained in \textbf{Step 2}, we know from Proposition \ref{pr: solva EDO paramatrique}, that for $\varepsilon<<1,~\kappa<<1$, there exists 
\begin{align*}
&\overline{S}=\overline{S}(\overline{\zeta},\varepsilon,\kappa,\eta,\gamma)=\overline{S}\Big(\overline{\zeta},\varepsilon,\kappa,|\overline{t}\big(\lambda, \Theta(f,g)\big)|+\eta,|h|,|\partial_x\overline{\phi}|,\overline{u}^\kappa(0),\big(\underline{u}_i^\kappa(\varepsilon)\big)_{i\in [\![1,I]\!]}\Big),~~\\
&\underline{S}=\underline{S}(\overline{\zeta},\varepsilon,\kappa,\eta,\gamma)=\overline{S}\Big(\overline{\zeta},\varepsilon,\kappa,|\underline{t}\big(\lambda, \Theta(f,g)\big)|+\eta,|h|,|\partial_x\underline{\phi}|,\underline{v}^\kappa(0),\big(\overline{v}_i^\kappa(\varepsilon)\big)_{i\in [\![1,I]\!]}\Big),\\  
\end{align*}
such that:
\begin{align}
&\label{eq absor derivee 1}
\partial_l\overline{\phi}(0,\overline{l}_\kappa)=\overline{S}(\overline{\zeta},\varepsilon,\kappa,\eta,\gamma)\ge \varepsilon I\overline{\zeta}\Big(|\overline{t}\big(\lambda, \Theta(f,g)\big)|+\eta+\frac{|b||\partial_x\overline{\phi}|+|h|}{\underline{\sigma}}\Big),\\ &\label{eq absor derivee 2}
\underline{\partial}_l\phi(0,\underline{l}_\kappa)=-\underline{S}(\overline{\zeta},\varepsilon,\kappa,\eta,\gamma)\leq - \varepsilon I\overline{\zeta}\Big(|\underline{t}\big(\lambda, \Theta(f,g)\big)|+\eta+\frac{|b||\partial_x\underline{\phi}|+|h|}{\underline{\sigma}}\Big).
\end{align}
Recall that $\overline{\zeta}$ is given in assumption $(\mathcal{H}-\bf R)$, and we have denoted for $w\in \mathcal{C}^2\big(\mathcal{N}_\varepsilon\times [0,\kappa]\big)$, $|\partial_xw|=\max_i\sup_{(x,l)}|\partial_xw_i(x,l)|$.
As an interpretation, roughly speaking \eqref{eq absor derivee 1} and \eqref{eq absor derivee 2} implies that the 'local time' time derivative of the test function $\big(\overline{\phi},\underline{\phi}\big)$ at $\bf 0$ absorb the error term induced by the Kirchhoff's speed of the Hamlitonians defined in \eqref{eq Kir speed}, scaled by $\varepsilon$.

Remark first that since $\overline{l}_\kappa \in [l_\star-\kappa,l_\star+\kappa]$ and $\overline{l}_\kappa \in [l_\star-\kappa,l_\star+\kappa]$, we have:
\begin{eqnarray}\label{eq cv point}
 \lim_{\kappa \searrow 0} \overline{l}_\kappa =l_\star,~~ \lim_{\kappa \searrow 0} \underline{l}_\kappa =l_\star
\end{eqnarray}
and from the expressions given in \eqref{eq  expression const scal} and the continuity of $u$ and $v$ we get:
\begin{eqnarray}\label{eq cv point bord}
&\nonumber \forall i\in[\![1,I]\!],~~\lim_{\varepsilon \searrow 0}\limsup_{\kappa \searrow 0} \overline{v}^\kappa_i(\varepsilon) =v(0,l_\star),~~\lim_{\varepsilon \searrow 0}\limsup_{\kappa \searrow 0} \underline{u}^\kappa_i(\varepsilon) =u(0,l_\star),\\
&\lim_{\kappa \searrow 0} \underline{v}^\kappa(0) =v(0,l_\star),~~\lim_{\kappa \searrow 0} \overline{u}^\kappa(0) =u(0,l_\star).
\end{eqnarray}
Using \eqref{eq: contrad visco en 0} and \eqref{eq fonction scaller}, we know that for all $i\in[\![1,I]\!]$:
$$u_i(\varepsilon,l_\star)-u(0,l_\star)-v_i(\varepsilon,l_\star)+v(0,l_\star)=f_i(\varepsilon,l_\star)-f(0,l_\star)-g_i(\varepsilon,l_\star)+g(0,l_\star)\ge 0.$$
Hence we deduce that for all $i\in[\![1,I]\!]$:
$$\lim_{\kappa \searrow 0}~~
\underline{u}^\kappa_i(\varepsilon)-\overline{u}^\kappa(0)-\overline{v}^\kappa_i(\varepsilon)+\underline{v}^\kappa(0)\ge 0,$$
and therefore:
\begin{align} \label{eq positivite u-v}
\exists \kappa_\varepsilon>0,~~\forall \kappa \leq \kappa_\varepsilon,~~
\underline{u}_i^\kappa(\varepsilon)-\overline{u}^\kappa(0)-\overline{v}^\kappa_i(\varepsilon)+\underline{v}^\kappa(0) \ge -\varepsilon^3.   
\end{align}
From the expressions of the ordinary differential equations \eqref{expression EDO sur sol}-\eqref{expression EDO sous sol} satisfied by that $\overline{\phi}$ and $\underline{\phi}$, we obtain:
\begin{align}\label{eq expression derivee en 0 fonctions tests}
\nonumber &\partial_x\overline{\phi}_i(0,\overline{l}_\kappa)=-\frac{\gamma}{\varepsilon}+\frac{1}{\varepsilon}\Big(\underline{u}^\kappa_i(\varepsilon)-\overline{u}^\kappa(0)-\int_0^\varepsilon\int_0^u\frac{\lambda \overline{\phi}_i(z,\overline{l}_\kappa)}{\overline{\sigma}}dzdu\\ &\nonumber -\int_0^\varepsilon\int_0^u\big[\overline{t}\big(\lambda, \Theta(f,g)\big)+\frac{|b||\partial_x\overline{\phi}_i(z,\overline{l}_\kappa)|+|h|}{\underline{\sigma}}+\eta\big]dzdu\Big),\\
&\nonumber \partial_x\underline{\phi}_i(0,\underline{l}_\kappa)=\frac{\gamma}{\varepsilon}+\frac{1}{\varepsilon}\Big(\overline{v}^\kappa_i(\varepsilon)-\underline{v}^\kappa(0)-\int_0^\varepsilon\int_0^u\frac{\lambda \underline{\phi}_i(z,\underline{l}_\kappa)}{\overline{\sigma}}dzdu\\
&\nonumber -\int_0^\varepsilon\int_0^u\big[\underline{t}\big(\lambda, \Theta(f,g)\big)-\frac{|b||\partial_x\underline{\phi}_i(z,\underline{l}_\kappa)|+|h|}{\underline{\sigma}}-\eta\big]dzdu\Big)\\
&\partial_l\overline{\phi}(0,\overline{l}_\kappa)=\overline{S}(\overline{\zeta},\varepsilon,\kappa,\eta,\gamma),~~~~
\partial_l\underline{\phi}(0,\underline{l}_\kappa)=-\underline{S}(\overline{\zeta},\varepsilon,\kappa,\eta,\gamma).
\end{align}
Now we have the necessary tools to conclude. Recall that $\vartheta\in \mathcal{O}$ is fixed. Write:
\begin{align*}
&\partial_l\overline{\phi}(0,\overline{l}_\kappa)+\sum_{i=1}^I\mathbb{S}_i(\overline{l}_\kappa,\vartheta)\partial_x\overline{\phi}_i(0,\overline{l}_\kappa)+h_0(\overline{l}_\kappa,\vartheta)-\\
&\Big(\partial_l\underline{\phi}(0,\underline{l}_\kappa)+\sum_{i=1}^I\mathbb{S}_i(\underline{l}_\kappa,\vartheta)\partial_x\underline{\phi}_i(0,\underline{l}_\kappa)+h_0(\underline{l}_\kappa,\vartheta)\Big)=\\
&\overline{S}(\overline{\zeta},\varepsilon,\kappa,\eta,\gamma)-\sum_{i=1}^I\mathbb{S}_i(\overline{l}_\kappa,\vartheta)\frac{1}{\varepsilon}\Big(\int_0^\varepsilon\int_0^u\big[\overline{t}\big(\lambda, \Theta(f,g)\big)+\\
&\frac{|b||\partial_x\overline{\phi}_i(z,\overline{l}_\kappa)|+|h|}{\underline{\sigma}}+\eta\big]dzdu\Big)+\underline{S}(\overline{\zeta},\varepsilon,\kappa,\eta,\gamma)-\\
&\sum_{i=1}^I\mathbb{S}_i(\underline{l}_\kappa,\vartheta)\frac{1}{\varepsilon}\Big(\int_0^\varepsilon\int_0^u\big[-\underline{t}\big(\lambda, \Theta(f,g)\big)+\frac{|b||\partial_x\underline{\phi}_i(z,\underline{l}_\kappa)|+|h|}{\underline{\sigma}}+\eta \big]dzdu \Big)-\\
&\frac{\gamma}{\varepsilon}\Big(\sum_{i=1}^I\mathbb{S}_i(\overline{l}_\kappa,\vartheta)+\sum_{i=1}^I\mathbb{S}_i(\underline{l}_\kappa,\vartheta)\Big)+\frac{1}{\varepsilon}\sum_{i=1}^I\mathbb{S}_i(\overline{l}_\kappa,\vartheta)\Big(\underline{u}^\kappa_i(\varepsilon)-\overline{u}^\kappa(0)-\overline{v}^\kappa_i(\varepsilon)+\underline{v}^\kappa(0)\Big)+\\
&\frac{1}{\varepsilon}\sum_{i=1}^I\Big(\mathbb{S}_i(\overline{l}_\kappa,\vartheta)-\mathbb{S}_i(\underline{l}_\kappa,\vartheta)\Big)\Big(\overline{v}^\kappa_i(\varepsilon)-\underline{v}^\kappa(0)\Big)+ h_0(\overline{l}_\kappa,\vartheta)-h_0(\underline{l}_\kappa,\vartheta)+\\
&\sum_{i=1}^I\mathbb{S}_i(\underline{l}_\kappa,\vartheta)\frac{1}{\varepsilon}\int_0^\varepsilon\int_0^u\frac{\lambda \underline{\phi}_i(z,\underline{l}_\kappa)}{\overline{\sigma}}dzdu-
\sum_{i=1}^I\mathbb{S}_i(\overline{l}_\kappa,\vartheta)\frac{1}{\varepsilon}\int_0^\varepsilon\int_0^u\frac{\lambda \overline{\phi}_i(z,\overline{l}_\kappa)}{\overline{\sigma}}dzdu.
\end{align*}
Now using \eqref{eq absor derivee 1}-\eqref{eq absor derivee 2}-\eqref{eq positivite u-v}, together with assumption $(\mathcal{H}-\bf R)$ we obtain that there exists $\exists \kappa_\varepsilon>0,~~\forall \kappa \leq \kappa_\varepsilon$:
\begin{align}\label{eq inegalite I1}
&\nonumber\partial_l\overline{\phi}(0,\overline{l}_\kappa)+\sum_{i=1}^I\mathbb{S}_i(\overline{l}_\kappa,\vartheta)\partial_x\partial_x\overline{\phi}_i(0,\overline{l}_\kappa)+h_0(\overline{l}_\kappa,\vartheta)-\\
&\nonumber\Big(\partial_l\underline{\phi}(0,\underline{l}_\kappa)+\sum_{i=1}^I\mathbb{S}_i(\overline{l}_\kappa,\vartheta)\partial_x\underline{\phi}_i(0,\underline{l}_\kappa)+h_0(\underline{l}_\kappa,\vartheta)\Big)\\
&\nonumber\ge -\frac{\gamma}{\varepsilon}2I\overline{\zeta}-\varepsilon^2I\overline{\zeta}-\frac{2|v|I\overline{\zeta}}{\varepsilon}|\overline{l}_\kappa-\underline{l}_\kappa|-|h||\overline{l}_\kappa-\underline{l}_\kappa|\\
&+\sum_{i=1}^I\mathbb{S}_i(\underline{l}_\kappa,\vartheta)\frac{1}{\varepsilon}\int_0^\varepsilon\int_0^u\frac{\lambda \underline{\phi}_i(z,\underline{l}_\kappa)}{\overline{\sigma}}dzdu-
\sum_{i=1}^I\mathbb{S}_i(\overline{l}_\kappa,\vartheta)\frac{1}{\varepsilon}\int_0^\varepsilon\int_0^u\frac{\lambda \overline{\phi}_i(z,\overline{l}_\kappa)}{\overline{\sigma}}dzdu. 
\end{align}
Taking now the infimum over all the $\vartheta \in \mathcal{O}$ in \eqref{eq inegalite I1}, we obtain with aid of \eqref{eq second jet u} and \eqref{eq second jet v}:
\begin{align}\label{eq inegalite I2}
\nonumber &0\ge -\frac{\gamma}{\varepsilon}2I\overline{\zeta}-\varepsilon^2I\overline{\zeta}-\frac{2|v|I\overline{\zeta}}{\varepsilon}|\overline{l}_\kappa-\underline{l}_\kappa|-|h||\overline{l}_\kappa-\underline{l}_\kappa|\\
&+\inf_{\vartheta \in \mathcal{O}}\Big[\sum_{i=1}^I\mathbb{S}_i(\underline{l}_\kappa,\vartheta)\frac{1}{\varepsilon}\int_0^\varepsilon\int_0^u\frac{\lambda \underline{\phi}_i(z,\underline{l}_\kappa)}{\overline{\sigma}}dzdu-
\sum_{i=1}^I\mathbb{S}_i(\overline{l}_\kappa,\vartheta)\frac{1}{\varepsilon}\int_0^\varepsilon\int_0^u\frac{\lambda \overline{\phi}_i(z,\overline{l}_\kappa)}{\overline{\sigma}}dzdu\Big]. 
\end{align}
For all $n>0$, there exists $\vartheta_n=\vartheta_n(\varepsilon,\kappa,\eta,\gamma)$ such that:
\begin{align*}
&\nonumber \inf_{\vartheta \in \mathcal{O}}\Big[\sum_{i=1}^I\mathbb{S}_i(\underline{l}_\kappa,\vartheta)\frac{1}{\varepsilon}\int_0^\varepsilon\int_0^u\frac{\lambda \underline{\phi}_i(z,\underline{l}_\kappa)}{\overline{\sigma}}dzdu-\\
&
\sum_{i=1}^I\mathbb{S}_i(\overline{l}_\kappa,\vartheta)\frac{1}{\varepsilon}\int_0^\varepsilon\int_0^u\frac{\lambda \overline{\phi}_i(z,\overline{l}_\kappa)}{\overline{\sigma}}dzdu\Big]+1/n\ge \\
&\Big[\sum_{i=1}^I\mathbb{S}_i(\underline{l}_\kappa,\vartheta_n)\frac{1}{\varepsilon}\int_0^\varepsilon\int_0^u\frac{\lambda \underline{\phi}_i(z,\underline{l}_\kappa)}{\overline{\sigma}}dzdu-
\sum_{i=1}^I\mathbb{S}_i(\overline{l}_\kappa,\vartheta_n)\frac{1}{\varepsilon}\int_0^\varepsilon\int_0^u\frac{\lambda \overline{\phi}_i(z,\overline{l}_\kappa)}{\overline{\sigma}}dzdu\Big].
\end{align*}
Dividing \eqref{eq inegalite I2} by $\varepsilon/2$ we get:
\begin{align}
&\nonumber 0\ge-\frac{\gamma}{\varepsilon^2}4I\overline{\zeta}-\varepsilon2I\overline{\zeta}-\frac{|v|4I\overline{\zeta}}{\varepsilon^2}|\overline{l}_\kappa-\underline{l}_\kappa|-\frac{2}{\varepsilon}|h||\overline{l}_\kappa-\underline{l}_\kappa|-\frac{2}{n\varepsilon}+\\
&\nonumber \Big[\sum_{i=1}^I\mathbb{S}_i(\underline{l}_\kappa,\vartheta_n)\frac{2}{\varepsilon^2}\int_0^\varepsilon\int_0^u\frac{\lambda \underline{\phi}_i(z,\underline{l}_\kappa)}{\overline{\sigma}}dzdu-
\sum_{i=1}^I\mathbb{S}_i(\overline{l}_\kappa,\vartheta_n)\frac{2}{\varepsilon^2}\int_0^\varepsilon\int_0^u\frac{\lambda \overline{\phi}_i(z,\overline{l}_\kappa)}{\overline{\sigma}}dzdu\Big]\ge\\
&\label{eq F1}-\frac{\gamma}{\varepsilon^2}4I\overline{\zeta}-\varepsilon2I\overline{\zeta}-\frac{|v|4I\overline{\zeta}}{\varepsilon^2}|\overline{l}_\kappa-\underline{l}_\kappa|-\frac{2}{\varepsilon}|h||\overline{l}_\kappa-\underline{l}_\kappa|-\frac{2}{n\varepsilon}+\\
&\label{eq F2}\frac{\lambda}{\overline{\sigma}}\underline{v}^\kappa(0)\sum_{i=1}^I\mathbb{S}_i(\underline{l}_\kappa,\vartheta_n)-\frac{\lambda}{\overline{\sigma}}\overline{u}^\kappa(0)
\sum_{i=1}^I\mathbb{S}_i(\overline{l}_\kappa,\vartheta_n)+\\
&\label{eq F3}\sum_{i=1}^I\mathbb{S}_i(\underline{l}_\kappa,\vartheta_n)\frac{2}{\varepsilon^2}\int_0^\varepsilon\int_0^u\frac{\lambda\big(\underline{\phi}_i(z,\underline{l}_\kappa)-\underline{v}^\kappa(0)\big)}{\overline{\sigma}}dzdu-\\
&\label{eq F4}
\sum_{i=1}^I\mathbb{S}_i(\overline{l}_\kappa,\vartheta_n)\frac{2}{\varepsilon^2}\int_0^\varepsilon\int_0^u\frac{\lambda\big(\overline{\phi}_i(z,\overline{l}_\kappa)-\overline{u}^\kappa(0)\big)}{\overline{\sigma}}dzdu
\end{align}
To conclude we have first in \eqref{eq F1}:
$$\underset{{\gamma \searrow 0,n \to +\infty}}{\lim_{\varepsilon \searrow 0,\kappa \searrow 0}}\Big[-\frac{\gamma}{\varepsilon^2}4I\overline{\zeta}-\varepsilon2I\overline{\zeta}-\frac{|v|4I\overline{\zeta}}{\varepsilon^2}|\overline{l}_\kappa-\underline{l}_\kappa|-\frac{2}{\varepsilon}|h||\overline{l}_\kappa-\underline{l}_\kappa|-\frac{2}{n\varepsilon}\Big]=0.$$
Secondly, from \eqref{eq positiv u v} and \eqref{eq cv point bord}, we obtain that in \eqref{eq F2} (because $\underline{v}^\kappa(0)\ge 0$, and $\overline{u}^\kappa(0)\leq 0$):
$$\limsup_{\kappa \searrow 0}\Big(\frac{\lambda}{\overline{\sigma}}\underline{v}^\kappa(0)\sum_{i=1}^I\mathbb{S}_i(\underline{l}_\kappa,\vartheta_n)-\frac{\lambda}{\overline{\sigma}}\overline{u}^\kappa(0)
\sum_{i=1}^I\mathbb{S}_i(\overline{l}_\kappa,\vartheta_n)\Big)\ge I\underline{\zeta}\frac{\lambda}{\overline{\sigma}}\big[v(0,l_\star)-u(0,l_\star)\big].$$
Finally, for the term \eqref{eq F3}, since the boundary condition of $\big(\overline{\phi}=(\varepsilon,\kappa,\eta,\gamma),\overline{\phi}=(\varepsilon,\kappa,\eta,\gamma)\big)$ satisfy the point (iii) of Proposition \ref{pr: solva EDO paramatrique}, namely (see  \eqref{eq cv point bord}):
$$\forall i\in[\![1,I]\!],~~ \lim_{\varepsilon \searrow 0}\limsup_{\kappa \searrow 0}|\underline{u}_i^\kappa(\varepsilon)-\overline{u}^\kappa(0)|=0,~~\lim_{\varepsilon \searrow 0}\limsup_{\kappa \searrow 0}|\overline{v}_i^\kappa(\varepsilon)-\underline{v}^\kappa(0|=0,$$
we have:
\begin{align*}
&\underset{\eta\searrow 0,\gamma \searrow 0}{\limsup_{\varepsilon \searrow 0,\kappa \searrow 0,}}\Big|\sum_{i=1}^I\mathbb{S}_i(\underline{l}_\kappa,\vartheta_n)\frac{2}{\varepsilon^2}\int_0^\varepsilon\int_0^u\frac{\lambda\big(\underline{\phi}^{\varepsilon,\kappa,\eta,\gamma}_i(z,\underline{l}_\kappa)-\underline{v}^\kappa(0)\big)}{\overline{\sigma}}dzdu\Big|\leq \\
&I\overline{\zeta}~~\underset{\eta\searrow 0,\gamma \searrow 0}{\limsup_{\varepsilon \searrow 0,\kappa \searrow 0,}}\max_i\underset{l\in [l_\star-\kappa,l_\star+\kappa]}{\sup}\Big|\frac{2}{\varepsilon^2}\int_0^\varepsilon\int_0^u\frac{\lambda\big(\underline{\phi}^{\varepsilon,\kappa,\eta,\gamma}_i(z,l)-\underline{v}^\kappa(0)\big)}{\overline{\sigma}}dzdu\Big|=0,
\end{align*}
and also:
\begin{align*}
&\underset{\eta\searrow 0,\gamma \searrow 0}{\limsup_{\varepsilon \searrow 0,\kappa \searrow 0,}}\big|\sum_{i=1}^I\mathbb{S}_i(\overline{l}_\kappa,\vartheta_n)\frac{2}{\varepsilon^2}\int_0^\varepsilon\int_0^u\frac{\lambda\big(\overline{\phi}^{\varepsilon,\kappa,\eta,\gamma}_i(z,\overline{l}_\kappa)-\overline{u}^\kappa(0)\big)}{\overline{\sigma}}dzdu\Big|=0.
\end{align*}
At the end we get:
\begin{align*}
&0\ge I\underline{\zeta}\frac{\lambda}{\overline{\sigma}}(v(0,l_\star)-u(0,l_\star)).    
\end{align*}
Therefore \eqref{eq fonction scaller} imply
$$0\ge g(0_,l_\star)-f(0,l_\star),$$
and that leads to a contradiction with \eqref{eq: contrad visco en 0}. We conclude that for all $\underline{\ell}\in (0,K)$, for all $(x,i)\in \mathcal{N}_R$ and for all $l\in [\underline{\ell},K]$:
$$f_i(x,l)\ge g_i(x,l).$$
We conclude the proof using the continuity of $f$ and $g$ with respect to variable $l$.
\end{proof}
\textbf{Proof of Theorem \ref{th: compa theorem l infiny}}:
\begin{proof}
 Let $u$ (and resp. $v$) be a super (resp. sub) continuous viscosity solution of \eqref{eq HJB bord infini}. Fix in the sequel $\alpha>0$ and $\underline{\ell}\in (0,+\infty)$. We argue by contradiction assuming that:
 $$\sup\Big\{~v_i(x,l)-u_i(x,l)-\alpha l^2,~~\big((x,i),l\big)\in \mathcal{N}_R\times[\underline{\ell},+\infty)~\Big\}>0.$$
 Using the growth condition satisfied by $u$ and $v$  with respect to the variable $l$ and the boundary conditions, it follows that the last supremum is necessary reached at a point:
 $$\big(x_\star,i_\star,l_\star\big)\in [0,R)\times[\![1,I]\!]\times [\underline{\ell},+\infty).$$
 Moreover from:
 $$v_{i_\star}(x_\star,l_\star)-u_{i_\star}(x_\star,l_\star)-\alpha l_\star^2>0,$$
It follows that necessary we have:
$$v_{i_\star}(x_\star,l_\star)-u_{i_\star}(x_\star,l_\star)>0.$$
We can argue then like in the last proof of Theorem \ref{th: compa theorem principale} to obtain a contradiction.
Therefore we obtain:
$$\forall \big((x,i),l\big)\in \mathcal{N}_R\times[\underline{\ell},+\infty),~~v_i(x,l)\leq u_i(x,l)+\alpha l^2.$$
Sending $\alpha \searrow 0$, and using the continuity of $u$ and $v$ with respect to the variable $l$, we obtain the required result.
\end{proof}

\end{document}